\documentclass[oneside,draft,11pt]{amsart}
\usepackage{amssymb, amscd}
\usepackage[all]{xy}

\theoremstyle{plain}
\newtheorem{theorem}{Theorem}[section]
\newtheorem{lemma}[theorem]{Lemma}
\newtheorem*{thm}{Theorem}
\newtheorem*{theorem*}{}
\newtheorem{proposition}[theorem]{Proposition}
\newtheorem{corollary}[theorem]{Corollary}
\newtheorem*{cor}{Corollary}

\theoremstyle{definition}
\newtheorem{conjecture}[theorem]{Conjecture}
\newtheorem{definition}[theorem]{Definition}
\newtheorem{problem}[theorem]{Problem}
\newtheorem{example}[theorem]{Example}

\theoremstyle{remark}
\newtheorem{remark}[theorem]{Remark}

\DeclareMathOperator{\Tr}{tr}
\DeclareMathOperator{\Id}{Id}

\DeclareMathOperator{\Tor}{Tor}
\DeclareMathOperator{\Hom}{Hom}
\DeclareMathOperator{\Res}{res}
\DeclareMathOperator{\Gal}{Gal}
\DeclareMathOperator{\dif}{d}

\DeclareMathOperator{\Image}{Im}

\DeclareMathOperator{\Sing}{{\mathcal S}}

\newcommand{\Sone}{{\mathbb S}^1}

\begin{document}

\title{On the Cohomology of Galois Groups Determined by Witt Rings}

\author{Alejandro Adem}
\address{Mathematics Department\\
         University of Wisconsin\\
         Madison, Wisconsin, 53706} 
\email{adem@math.wisc.edu}
\thanks{The first author was partially supported by an NSF grant and the MPIM-Bonn.}

\author{Dikran B. Karagueuzian}
\address{Mathematics Department\\
         University of Wisconsin\\
         Madison, Wisconsin, 53706} 
\email{dikran@math.wisc.edu}
\thanks{The second author was 
partially supported by an NSF postdoctoral fellowship, CRM-Barcelona,
and MPIM-Bonn.}

\author{Jan Minac}
\address{Mathematics Department \\
         University of Western Ontario \\
         London, Ontario, Canada N6A 5B7}
\email{minac@uwo.ca}
\thanks{The third author was partially supported by the NSERC and the
special Dean of Science fund at UWO. He would like to thank Wenfeng
Gao for his influence as a collaborator during the past few years.
Gao's original ideas were an important inspiration for these
investigations on \textit{W}-groups.}


\begin{abstract}

Let $F$ denote a field of characteristic different from two.
In this paper
we describe the mod 2 cohomology of a Galois group ${\mathcal G}_F$ (called
the \textit{W}-group of $F$) which is known to essentially
characterize the Witt ring
$W\!F$
of anisotropic quadratic modules over $F$.
We show that $H^*({\mathcal G}_F, {\mathbb F}_2)$ contains the mod 2 Galois cohomology of
$F$ and that its structure will reflect important properties of the field.  
We construct a space $X_F$ endowed with an action of an elementary abelian
group $E$ such that the computation of the cohomology of ${\mathcal G}_F$
reduces to calculating the equivariant cohomology
$H^*_E(X_F,{\mathbb F}_2)$. For the case of a field which is not
formally real this amounts to computing 
the cohomology of an explicit Euclidean space form, an object which is
interesting in its own right.
We provide a number of
examples and a substantial combinatorial computation for the cohomology
of the universal \textit{W}-groups.
\end{abstract}

\maketitle

\begin{section}{Introduction}
\label{s:introduction}

Although there has been substantial recent activity in the area of
finite group cohomology, it has hardly involved any interactions
with Galois theory. This seems rather surprising in light of the
recent major developments in field theory, such as Voevodsky's
proof of the Milnor Conjecture \cite{V}. Moreover, cohomological
methods in field theory use comparatively little from 
available techniques in the
cohomology of finite groups. An
explanation for this is the fact that interesting Galois groups are
usually far too complicated to be analyzed with basic cohomological
methods. An obvious compromise would be to identify quotient 
groups which still retain substantial field-theoretic information
but which are accessible via methods in the cohomology of groups.

In \cite{MS1}, Minac and Spira introduced a relatively simple
Galois group ${\mathcal G}_F$
associated to a field $F$ of characteristic different
from $2$, known as the \textit{W}-group of $F$. They showed that, up to a minor
technical condition, this group will characterize the Witt ring
$W\! F$ of anisotropic quadratic modules over $F$
(see \cite{M} for background). Hence these groups would seem to be ideal
candidates for a fruitful cohomological analysis, where ideas and methods
from Galois theory could be successfully combined with well-established
techniques from the cohomology of groups.
In this paper we show that this is indeed the case: after showing
that ${\mathcal G}_F$ can be intrinsically defined as a central
extension of an elementary abelian $2$-group by an elementary
abelian $2$-group using
field-theoretic data, we prove that its mod $2$ cohomology
not only contains the Galois cohomology but that its qualitative
structure reflects properties of the field. We also provide a
topological method for computing $H^*({\mathcal G}_F,{\mathbb F}_2)$,
reducing matters in many instances to determining the cohomology of 
a compact Euclidean space form, hence allowing us to apply
methods from homotopy theory, group actions and
combinatorics in our analysis. 

Before stating our main results, we recall some
basic definitions. Throughout this paper, $F$ will denote a field of
characteristic \emph{different from two}. This point is emphasized
here as it will be an unstated hypothesis in most of the results in
the rest of this paper. We denote the quadratic closure of $F$ by
$F_q$ and the associated Galois group by $G_F$.  The \textit{W}-group of $F$
is defined (see \cite{MS1}) as the quotient group 
\[
{\mathcal G}_F = G_F/G_F^4[G_F^2,G_F].
\]
A fairly simple description of the intermediate extension 
$F\subset F^{(3)}\subset F_q$ such that ${\mathcal
G}_F=\Gal(F^{(3)}/F)$ is given in \S \ref{s:w-groups}.
Let $\Phi ({\mathcal G}_F)
\subset {\mathcal G}_F$ denote its Frattini subgroup, then it is an elementary
abelian $2$-group and ${\mathcal G}_F$ can be expressed as a central extension
\[
1\to \Phi ({\mathcal G}_F)\to {\mathcal G}_F\to E\to 1
\]
where $E\cong (\dot F/\dot F^2)^*$, the Pontrjagin dual of the
mod $2$ vector space $\dot F/\dot F^2$.
The starting point for our work is
a field-theoretic characterization of a condition on groups 
(known as the 2C property, see \cite{AK}) which
has substantial cohomological implications:
\begin{thm} {\rm \textbf{(\ref{t:2C})}}
If $|{\mathcal G}_F|>2$, then
$-1\in F$ is a sum of squares ($F$ is not formally real)
if and only if every element of order $2$ in
${\mathcal G}_F$ is central.
\end{thm}

In terms of the defining extension this implies that $\Phi ({\mathcal G}_F)$
is the unique maximal elementary abelian subgroup
of ${\mathcal G}_F$; in the formally
real case we show that $\Phi ({\mathcal G}_F)$ is an index two subgroup
of any 
maximal elementary abelian subgroup.

Given this simple expression for ${\mathcal G}_F$ as a central
extension, the next logical step is to identify it
using cohomological and field-theoretic data. Recall that if
\[
l\colon \dot F/\dot F^2\to K_1F/2K_1F
\]
is the canonical isomorphism between $\dot F/\dot F^2$ written multiplicatively
and additively respectively, then Milnor \textit{K}-theory mod 2 (see \cite{M}) can be
expressed as 
${\mathbb F}_2 [x_i \mid i\in \Omega ]/I_F$, where $\{x_i\mid i\in\Omega\}$
are one dimensional polynomial generators which constitute a basis for 
$K_1F/2K_1F$ 
and $I_F$ is the ideal generated by the
quadratic polynomials corresponding to $l(a)l(1-a)$, for $a\in \dot
F/\dot F^2$, 
$a\ne 1$. Let ${\mathcal B}_F$ denote the subspace of $H^2(E, {\mathbb
F}_2)$ spanned by these 
polynomials. From the five term exact sequence associated to the 
defining extension for ${\mathcal G}_F$,
we obtain an injective map $\delta \colon H^1(\Phi ({\mathcal G}_F),{\mathbb
F}_2) \to H^2(E,{\mathbb F}_2)$. 

\begin{thm} {\rm \textbf{(\ref{t:central})}} 
The image of $\delta$ is the subspace ${\mathcal B}_F\subset
H^2(E,{\mathbb F}_2)$, and  
${\mathcal G}_F$ is the uniquely
determined central extension
associated to this subspace.
\end{thm}

One can think of ${\mathcal G}_F$ as the ``minimal'' extension where
the relations in \textit{K}-theory are satisfied.
As a consequence of this and the main result in \cite{MS1}
we obtain 

\begin{cor}{\rm \textbf{(\ref{t:corollary})}}
If $F$ and $L$ are fields as above
then ${\mathcal G}_F\cong {\mathcal G}_L$
if and only if $K_1F/2K_1F\cong K_1L/2K_1L$. Hence if
$WF\cong WL$, $K_1F/2K_1F\cong K_1L/2K_1L$, and the converse
is true provided $F$ and $L$ have the same level if 
$\langle 1,1\rangle_F$ is universal.
\end{cor}

Next we use the Milnor
Conjecture to prove 

\begin{thm} {\rm \textbf{(\ref{2.3})}} 
Let ${\mathcal R}\subset 
H^*({\mathcal G}_F,{\mathbb F}_2)$ denote the subring generated by
one dimensional classes, then 
${\mathcal R}\cong H^*(F, {\mathbb F}_2)$, the mod 2 Galois Cohomology
of $F$.
\end{thm}

An important point to make is that 
$H^*({\mathcal G}_F, {\mathbb F}_2)$ can contain substantially 
more cohomology than $H^*(F, {\mathbb F}_2)$ alone, and in fact its qualitative structure will
reflect properties of $F$. 
This is best illustrated when $|\dot F/\dot F^2|<\infty$,
an assumption which we make for the rest of this introduction.
If $|\dot F/\dot F^2|=2^n$, then $E$ and
$\Phi ({\mathcal G}_F)$ are finite elementary abelian groups
of ranks $n$ and $r$ respectively, where 
$r=\binom{n+1}{2}-\hbox{dim}~H^2(\hbox{Gal}(F_q/F), {\mathbb F}_2)$.
Using methods from the cohomology of finite groups,
we prove 

\begin{thm} {\rm \textbf{(\ref{t:cohen-macaulay}, \ref{t:not-quite-CM})}} 
There exist polynomial classes 
$$\zeta_1,\dots ,\zeta_r\in H^2({\mathcal G}_F, {\mathbb F}_2)$$
which form a regular sequence. If the field $F$ is  not formally real, then
$H^*({\mathcal G}_F, {\mathbb F}_2)$ is free and finitely generated over 
the polynomial subalgebra ${\mathbb F}_2 [\zeta_1,\ldots ,\zeta_r]$ 
(in particular Cohen-Macaulay). In the formally
real case, $H^*({\mathcal G}_F, {\mathbb F}_2)$ has depth equal to $r$ or $r+1$
and Krull dimension equal
to $r+1$.
\end{thm}

Using the results in \cite{AK}, we show

\begin{thm} {\rm \textbf{(\ref{p:essential-exterior})}} 
Let $F$ be a field which is not formally real. Then there exist
non-zero classes 
$x\in H^*({\mathcal G}_F, {\mathbb F}_2)$ 
which restrict trivially on all proper subgroups of
${\mathcal G}_F$. Any such class must be exterior, i.e.\ $x^2=0$.
\end{thm}

The undetectable classes are also called \textsl{essential} cohomology
classes; in fact they constitute an ideal 
${\mathcal E}\subset H^*({\mathcal G}_F,{\mathbb F}_2)$ such that any two
elements in it multiply trivially. Note that the result above
indicates that the usual detection methods for computing
cohomology will not work for these \textit{W}-groups. 
In contrast, for formally
real fields detection can occur and in fact the element 
$[-1]\in H^1({\mathcal G}_F, {\mathbb F}_2)$
plays a key role. Recall that a field is said to be pythagorean
if $F^2+F^2=F^2$. Our main result for formally real
fields is the following:

\begin{thm} {\rm \textbf{(\ref{2.4})}} 
Let $F$ be a formally real field, then the following conditions are
equivalent:
\begin{enumerate}
\item $[-1]\in H^1({\mathcal G}_F, {\mathbb F}_2)$ is not a zero divisor.
\item $F$ is pythagorean
and $H^*({\mathcal G}_F, {\mathbb F}_2)$ is Cohen-Macaulay.
\item $F$ is pythagorean and 
$H^*({\mathcal G}_F, {\mathbb F}_2)$ is detected on
its elementary abelian subgroups.
\end{enumerate}
\end{thm}

To calculate the quotient algebra
$H^*({\mathcal G}_F, {\mathbb F}_2)/(\zeta_1,\dots ,\zeta_r)$,
we construct
what we call a \textsl{topological model}, its main properties are summarized
in

\begin{thm} {\rm \textbf{(\ref{t:X_F})}} 
Given a field $F$ with $|\dot F/\dot F^2|=2^n$, there exists an action
of $E_n\cong ({\mathbb Z}/2)^n$ on $X_F\cong (\Sone)^r$ with the
following properties: 
\begin{enumerate}
\item $E_n$ only has cyclic isotropy subgroups,
\item the action is free if and only if $F$ is  not formally real,
\item $H^*({\mathcal G}_F, {\mathbb F}_2)/(\zeta_1,\dots
,\zeta_r)\cong H^*_{E_n}(X_F; {\mathbb F}_2).$
\end{enumerate}
\end{thm}

The term appearing in (3) is the \textsl{equivariant cohomology}, i.e.
the mod 2 cohomology of the \textsl{Borel construction}
$X_F\times_{E_n}EE_n$. In case (2) we of course obtain the cohomology
of $X_F/E_n$, a compact 
$r$-dimensional Euclidean space form 
such that the subring generated 
by one dimensional classes is isomorphic to $H^*(F,{\mathbb F}_2)$.
The fact that the $\zeta_i$ form a regular sequence implies \cite[10.3.4]{E}
that if $p_F(t)$ is the Poincar\'e series for the cohomology 
of ${\mathcal G}_F$, and $q_F(t)$ the one for the equivariant cohomology
$H^*_{E_n}(X_F, {\mathbb F}_2)$, then $p_F(t)=q_F(t)/(1-t^2)^r$.  If
$F$ is not formally real,
$q_F(t)$ will describe a basis for $H^*({\mathcal G}_F, {\mathbb F}_2)$
as a module over ${\mathbb F}_2 [\zeta_1 ,\dots ,\zeta_r]$.
Also we should point
out that the geometry of the action reflects the field theory in other
ways besides condition (2). For example, if $F$ is pythagorean, one 
can find a hyperplane $H\subset E_n$ which acts freely on $X_F$; this
will correspond to the index two subgroup ${\mathcal G}_{F(\sqrt{-1})}
\subset {\mathcal G}_F$.

We apply the results above to examples of interest in
field theory. For example, if $F={\mathbb Q}_2$ is the field of
$2$-adic numbers, we obtain a compact
$5$-dimensional manifold with Poincar\'e series equal to
$1+3t+6t^2+6t^3+3t^4+t^5$.
In the case of superpythagorean fields with $|\dot F/\dot F^2|<\infty$,
(such as $F_n
={\mathbb R}((t_1))((t_2))\dots ((t_{n-1}))$, the field of iterated power
series over ${\mathbb R}$), we obtain a complete description of the
cohomology of ${\mathcal G}_{F_n}$, and from there the Galois
cohomology of any superpythagorean field.

Given a \textit{W}-group ${\mathcal G}_F$ 
with $|\dot F/\dot F^2|=2^n$, it can be expressed as
a quotient of a unique
``universal'' \textit{W}-group $W(n)$ generated by $n$ elements.
In terms of the extension data this
is the group corresponding to the entire vector space $H^2(E,{\mathbb F}_2)$.
Given their defining properties, these groups will have the most
interesting and complicated cohomology. In particular $W(n)$
has the 2C property, hence its cohomology is ``undetectable''.
To compute $H^*(W(n), {\mathbb F}_2)$, we replace its given topological model
(an orbit space which is interesting in its own right)
by one which has much more plentiful
\textit{rational} cohomology.
Using an Eilenberg-Moore spectral sequence associated to 
this model, we obtain the following substantial combinatorial computation:

\begin{thm} {\rm \textbf{(\ref{t:combinatorics})}} 
If ${\mathcal G}_F=W(n)$ and $q_F(t)=1+a_1t+\dots +a_rt^r$, then
\[
a_i\ge \sum_{p+q=i} \sum_{Y_{\lambda}}\prod_{(s,t)\in Y_{\lambda}}
\frac {n+t-s}{h(s,t)}
\]
where $Y_{\lambda}$ ranges over all symmetric, $p+2q$-box,
$p$-hook Young diagrams, and $h(s,t)$ denotes the hooklength
of the box $(s,t)$.
\end{thm}

We can verify that this theorem
gives an equality for $i\le 3$, whence
we obtain
\[
a_1=n, \quad a_2=\frac{n(n+1)(n-1)}{3}, \quad a_3=\frac{n(n^2-1)(3n-4)(n+3)}{60}.
\]
It is worthwhile to note that our method uses rational techniques
to produce classes in the cohomology of a finite $2$-group, which
seems to be a somewhat novel approach.
Determining whether or not we have an equality for all
coefficients $a_1,\dots ,a_r$ is an interesting
problem, equivalent to showing that the homology of an integral Koszul
complex is
2-torsion free (see \S 6 for details).

The subsequent sections of this paper are organized as follows:
in \S 2 we provide the background
on \textit{W}-groups involving Galois theory; in \S 3 we discuss the
basic cohomological structure of the ${\mathcal G}_F$; in \S 4 we 
introduce our topological models together with examples, including
a discussion of the possible low-dimensional euclidean space forms
which can occur; in \S 5 we analyze the situation for formally real
fields, discussing at length the general pythagorean and
superpythagorean case; in \S 6 we provide the computation
for universal \textit{W}-groups; and finally in \S 7 we discuss a plausible 
general calculation for the cohomology of
\textit{W}-groups.

Throughout this paper coefficents will be assumed in the
field ${\mathbb F}_2$ with two elements unless stated
otherwise, hence they are suppressed from now on.
The results in this paper have been announced in \cite{AKM}.
We are grateful to J. Carlson, F. Cohen, D. Kotschick 
and J.-P. Serre for their useful
comments on aspects of this work.
\end{section}

\begin{section}{Preliminaries on \textit{W}-groups}
\label{s:w-groups}
In this section we will provide preliminary information 
on \textit{W}-groups, our basic reference is \cite{MS1}. 
Assume that $F$ is a field of characteristic
different from 2.
Let $F_q$ denote the
\textsl{quadratic closure} of $F$, and denote by
$G_F$ the Galois group of this extension over $F$. We will
begin by defining the \textit{W}-group associated to $F$.

\begin{definition}
Let $G_F^4[G_F^2,G_F]\subset G_F$ denote
the closure of the subgroup generated by fourth powers and
commutators of squares with arbitrary elements.
The \textit{W}-group of $F$ is the quotient group
\[
{\mathcal G}_{F}=G_F/G_F^4[G_F^2,G_F].
\]
\end{definition}

The group ${\mathcal G}_F$ is a pro-2-group, and it
is finite 
if and only if $|\dot F/\dot F^2|<\infty$.
Let $F^{(2)}=F(\sqrt{a} \mid a\in\dot F)$, ${\mathcal E}=
\{ y\in F^{(2)}\mid F^{(2)}(\sqrt{y})/F \quad \text{is Galois}\}$,
and $F^{(3)}=F^{(2)}(\sqrt{y} \mid y\in {\mathcal E})$. Then we have
a sequence of field extensions $F\subset F^{(2)}\subset F^{(3)}
\subset F_q$  where 
$\Gal(F^{(3)}/F)={\mathcal G}_{F}$,
$\Gal(F^{(2)}/F)\cong G_F/G_F^2\cong E$, an elementary abelian 2-group 
such that
$E\cong (\dot F/\dot F^2)^*$, the Pontrjagin dual, and $G_F^2$ is the Frattini
subgroup of $G_F$. 
Now if we let $\Phi ({\mathcal G}_F)\subset {\mathcal G}_F$
denote its Frattini subgroup, then it can be identified with
the subgroup $\Gal(F^{(3)}/F^{(2)})$, and we have a diagram of
extensions

\[
\xymatrix{
 &  1 \ar[d] & 1 \ar[d] \\
& G_F^4[G_F^2,G_F] \ar@{=}[r] \ar[d] & G_F^4[G_F^2,G_F] \ar[d]	\\
1 \ar[r] & G_F^2 \ar[r] \ar[d] & G_F \ar[r] \ar[d]	& E	\ar[r]
\ar@{=}[d] &	1 \\ 
1 \ar[r] & \Phi ({\mathcal G}_F) \ar[r] \ar[d] & {\mathcal G}_F \ar[r]
\ar[d] & E \ar[r] & 1\\ 
& 1 & 1 \\
}
\]
where the subgroup $\Phi ({\mathcal G}_F)$ is elementary abelian as well
as central.

Let $\mathcal C$ denote the class of finite $2$-groups
$H$ such that $H^4[H^2,H]=\{ 1\}$. If $\{e_i\mid i\in I\}$ is a basis for
the vector space $\dot F /\dot F^2$, then it is
not hard to show that ${\mathcal G}_{F}$ is a pro-$\mathcal C$-group
with a minimal set of generators of cardinality $|I|$.
The \textsl{universal} \textit{W}-group on $I$, $W(I)$, can be defined as
the unique pro-$\mathcal C$-group satisfying the following condition:
for any pro-$\mathcal C$-group $H$ and any (set) map $f \colon I \to
H$, there is a unique extension of $f$ to a pro-$\mathcal
C$-group-homomorphism $\tilde{f} \colon W(I) \to H$. 
 From \cite{MS1}, we know that
every element of order $2$ in $W(I)$ is central (i.e.\ $W(I)$
satisfies the 2C condition, see \cite{AK}) and that
its Frattini subgroup $\Phi (W(I))$ is a maximal elementary abelian
subgroup. Now if ${\mathcal G}_{F}$ is 
any \textit{W}-group, then it fits into an extension
$$1\to V\to W(I)\to {\mathcal G}_{F}\to 1$$
where $I$ is a set in one to one correspondence
with a minimal set of generators for
$\dot F/\dot F^2$ and $V\subseteq 
\Phi (W(I))$. In the special case when $|\dot F/\dot F^2|=2^n$,
$\Phi ({\mathcal G}_F)$ is an elementary abelian group
of rank equal to $r=n + {\binom{n}{2}} - \dim H^2(G_F)$

Recall that $F$ is said to be \textsl{formally real} if $-1$
is not a sum of squares in $F$. The following theorem
provides a characterization of \textit{W}-groups for 
fields which are not formally real.

\begin{theorem}
\label{t:2C}
Suppose that $|{\mathcal G}_F|>2$; then the field $F$
is not formally real if and only if each element of order
$2$ in ${\mathcal G}_F$ is central.
\end{theorem}

\begin{proof}
Assume first that $F$ is not formally
real. Let $\sigma$ be any element in
${\mathcal G}_F$ of order $2$. If $\sigma$ belongs to $\Phi ({\mathcal G}_F)$,
then it is central, as all the elements in this subgroup
are central. Therefore we may assume that $\sigma\notin\Phi ({\mathcal G}_F)$.
However, in \cite{MS2} it was shown that
\[
P_{\sigma}=\{f\in F \mid \sqrt{f}^{\sigma}=
\sqrt{f}\}
\]
is an ordering in the field $F$. Hence $F$ is a formally 
real field, a contradiction. Therefore if $F$ is not formally
real, condition 2C must hold.

Assume now that $F$ is formally real.
In \cite{MS2} it was shown that the correspondence
$\sigma\leftrightarrow P_{\sigma}$, $\sigma\in {\mathcal G}_F-\Phi ({\mathcal G}_F)$
where $\sigma\in {\mathcal G}_F$ is an involution, is a one-to-one
correspondence between classes of involutions not contained
in $\Phi ({\mathcal G}_F)$ and orderings in $F$. Hence we may assume the existence
of an involution $\sigma\in {\mathcal G}_F-\Phi ({\mathcal G}_F)$. 
Due to the fact that
$|{\mathcal G}_F|\ge 4$, there exists an element $\tau\in
{\mathcal G}_F-\{1,\sigma\}$. Since $\sigma$ is an involution, it
cannot generate the whole group ${\mathcal G}_F$. Therefore
$\tau$ can be chosen in ${\mathcal G}_F-\Phi ({\mathcal G}_F)$ and moreover
$\sigma$ and $\tau$ are linearly independent in the
vector space ${\mathcal G}_F/\Phi ({\mathcal G}_F)$.

As $F$ is formally real, $-1$ is not a square in $\dot F$. Moreover,
in \cite[p.521]{MS2}, it is proved that $(\sqrt{-1})^{\sigma}
=-\sqrt{-1}$. As $\sigma, \tau$ are linearly independent
modulo $\Phi ({\mathcal G}_F)$, we see that there exists an
element $a\in \dot F$ such that $(\sqrt{a})^{\sigma}=\sqrt{a}$
but $(\sqrt{a})^{\tau}=-\sqrt{a}$. Set $L=F(\root 4 \of a ,
\sqrt{-1})$; then $L/F$ is a Galois extension 
with Galois group $D_4$, the dihedral group of order eight.
 From the information we have it is easy to check that
the images of $\tau$ and $\sigma$ generate $D_4$ and in particular
they \emph{do not} commute. 
This implies that $\sigma$ is not central and we conclude that 
the 2C condition does not hold.
\end{proof}

\begin{remark}
In the paper
\cite{MMS} similar dihedral tricks are used to exclude certain
subgroups of ${\mathcal G}_F$. For later use we define the
\textsl{level} of a field $F$ as $s(F)$, the minimum $m$ such that 
there exist $f_1,\dots , f_m\in F$ with $-1=f_1^2 +\dots +f_m^2$.
Then $F$ is not formally real if and only if $s(F)<\infty$ and
Pfister has proved that $s(F)\in \{1,2,4,8,16,\dots ,\infty\}$.
\end{remark}

As a consequence of the theorem above, we see that if $F$ is not
formally real the Frattini subgroup $\Phi ({\mathcal G}_F)$ 
is the unique maximal elementary abelian subgroup, which in 
particular is central. In the formally real case we have

\begin{theorem}
\label{t:non2C}
If $F$ is formally real then $\Phi ({\mathcal G}_F)\subset
{\mathcal G}_F$ is a central elementary abelian subgroup of
index two in any maximal elementary abelian subgroup
in ${\mathcal G}_F$.
\end{theorem}

\begin{proof}
We see from the proof of \ref{t:2C} that any pair of involutions
in ${\mathcal G}_F-\Phi ({\mathcal G}_F)$ which are linearly independent
in ${\mathcal G}_F/\Phi ({\mathcal G}_F)$ cannot commute. On the other
hand as the 2C condition cannot hold we conclude that
the Frattini subgroup must be a central, index two subgroup of any
maximal elementary abelian subgroup of ${\mathcal G}_F$.
\end{proof}

As we shall see in the next section, the group theoretic results
in this section have important cohomological consequences. 
In particular
we will make use of the recent
cohomological characterization of finite $2$-groups with 
the 2C property provided in \cite{AK}.
\end{section}

\begin{section}{Group Cohomology and Extensions}
\label{s:cohomology}

We will require some basic facts from the cohomology of
finite groups, all of which can be found in \cite{AM} or \cite{E}.
First we have 

\begin{lemma} 
\label{l:varieties}
Let $G$ denote a finite $2$-group
satisfying the 2C condition with $E\subseteq G$ the elementary 
abelian subgroup of maximal rank $k$. Let ${\mathcal P}
\subset H^*(G)$ be a polynomial subalgebra such that $H^*(E)$
is a finitely generated module over $res^G_E({\mathcal P})$.
Then $H^*(G)$ is a free and finitely generated module over
${\mathcal P}$.
\end{lemma}

\begin{proof} Let ${\mathcal P}= {\mathbb F}_2[\zeta_1,\dots ,\zeta_n]$.
By a theorem of Duflot \cite{D}, we know that $H^*(G)$ is Cohen-Macaulay
and by a standard result in commutative algebra 
we know that under that condition the cohomology will be a free module
over any polynomial subring over which it is finitely generated. (For
a proof of this result see \cite[p.IV-20, Thm.2]{Serre-loc}; for a
historical discussion see \cite{EGH}.) Hence 
we only need to prove that $H^*(G)$ is a finitely generated $\mathcal
P$-module. To prove this we will use the more geometric language 
of cohomological varieties (see \cite{E} for background).

Let $V_G(\zeta_i)$ denote the homogeneous hypersurface in $V_G$,
(the maximal ideal spectrum for $H^*(G)$) defined by $\zeta_i$. Then
$H^*(G)$ will be finitely generated over $\mathcal P$ if and only if
$V_G(\zeta_1)\cap\dots
\cap V_G(\zeta_n)=\{0\}$. If we represent the class $\zeta_i$ by an
epimorphism $\Omega^{n_i}({\mathbb F}_2)\to {\mathbb F}_2$ 
with kernel $L_{\zeta_i}$, then 
we know that $V_G(\zeta_i)=V_G(L_{\zeta_i})$, the variety associated
to the annihilator of $Ext^*_{{\mathbb F}_2G}(L_{\zeta_i},L_{\zeta_i})$. Moreover 
using basic properties of these varieties, we have that
$V_G(\zeta_1)\cap\dots \cap V_G(\zeta_n)=V_G(L_{\zeta_1}\otimes\dots
\otimes L_{\zeta_i})$. Now the cohomological variety of a module will be
$0$ if and only if the module is projective, hence what we need to
prove is that the module
$L_{\zeta_1}\otimes\dots\otimes L_{\zeta_n}$ is projective. However
by Chouinard's Theorem (see \cite{E}) we know that it is enough to 
check this by restricting to maximal elementary abelian subgroups; in this case
$E$ is the only such group and projectivity follows from our
hypothesis, as $\Res_E^G({\mathcal P})\subseteq H^*(E)$ is a
polynomial subalgebra  
over which it is finitely generated (note that by Quillen's detection
theorem, the kernel of $\Res^G_E$ is nilpotent, hence $\mathcal P$ embeds in
$H^*(E)$ under this map).
Hence we conclude that $\zeta_1,\dots ,\zeta_n$ form
a ``homogeneous system of parameters'' and so $H^*(G)$ is free and finitely
generated as a module over $\mathcal P$. 
\end{proof}

Next we give a slight generalization of this lemma which is also
well-known to the experts (the idea is in \cite{D}),
but which is not available in the literature in precisely the form we
want.  Since the proof of the generalization is not so
conceptual as that of the lemma above, we confine ourselves to
indicating how an existing proof may be modified to 
yield the result we want.  

\begin{lemma}
\label{l:regular-sequence}
Let $E \subset G$ be an central elementary abelian $p$-subgroup and
let $\zeta_1, \ldots, \zeta_n$ be a sequence of elements in $H^*(G)$
which restrict to a regular sequence in $H^*(E)$. Then $\zeta_1,
\ldots, \zeta_n$ is a regular sequence in $H^*(G)$.
\end{lemma}

\begin{proof}
We will not need the case $p \neq 2$, so we ignore it. In
\cite[1.1]{BH}, this theorem is proved in the special case 
where there exists a basis $u_1, \ldots, u_n$ for $H^1(E)$ and
numbers $l_i$ such that $\Res_E(\zeta_i) = u_i^{l_i}$.  A brief
examination of the proof shows that the only property of these
elements $u_i^{l_i}$ used is that they form a regular sequence in
$H^*(E)$. (The reader who is not inclined to examine the proof in
\cite{BH} may be excused from doing so on the grounds that we will
actually be able to take our restrictions to have this special form.)
\end{proof}

Our next lemma is about group extensions.

\begin{lemma} 
\label{l:extension}
Let 
\[
1\to V\to G\to W\to 1
\]
be a central extension of elementary abelian $2$-groups where
$V=\Phi (G)$, the Frattini subgroup.
Then, the differential $\delta\colon H^1(V)\to H^2(W)$ 
in the five term exact sequence for the extension
above is a monomorphism, and if $G$ is finite,
the isomorphism class of the extension is determined
by the 
subspace $\delta (H^1(V))\subset H^2(W).$
\end{lemma}

\begin{proof} The fact that $\delta$ is injective follows directly
from the five-term exact sequence associated to the extension.
Now a 
well-known fact in group cohomology
is that the isomorphism class of a finite
extension as above is determined
by an \textsl{extension class} in $H^2(W, V)$. We will show that determining
this isomorphism class is equivalent to identifying the subspace
mentioned above. To obtain the extension class we consider the five term exact
sequence for the extension with coefficients in $V$; we have:
\[
0\to H^1(W,V)\to H^1(G, V)\to H^1(V,V)\buildrel {\delta}\over
\longrightarrow H^2(W, V)\to H^2(G, V).
\]
Our extension class is $\delta (\Id)$ (see \cite[p.~207]{HS}); we can
decompose it by using a basis for $H^1(V)$. 
This gives rise to a basis for
the subspace $\Image\delta$ and hence determines it uniquely. Conversely
up to a change of basis, the subspace $\Image\delta\subset H^2(W)$
determines the map 
$\delta\colon H^1(V)\to H^2(W)$ and hence $\delta (Id)$.
\end{proof}

Following standard usage, 
a basis for $im~\delta$ will be called a collection
of defining $k$-invariants for the group extension.

\begin{remark}
This result easily extends to pro-finite groups such
as our $W$-groups (see \cite{Ribes}, page 100). Hence
${\mathcal G}_F$ is uniquely determined by the extension data.
\end{remark}

\begin{remark} 
Consider the universal group on $n$ generators, $W(n)$, described as a
central extension
\[
1\to \Phi (W(n))\to W(n)\to E_n\to 1
\]
where $E_n\cong ({\mathbb Z}/2)^n$, and
$\Phi (W(n))\cong ({\mathbb Z}/2)^{n+\binom{n}{2}}$.
Then $W(n)$ is the central extension
associated to the \emph{entire} vector space $H^2(E_n)$.
\end{remark}

\begin{remark}
The mod 2 cohomology of a finite group is known to be Noetherian (see \cite{E}).
On the other hand it is elementary to verify that if $G$ is a pro-2-group
such that $G/\Phi (G)$ is infinite, then $H^*(G)$ cannot be Noetherian;
hence $H^*({\mathcal G}_F)$ is Noetherian if and only if 
$|\dot F/\dot F^2|<\infty$.
\end{remark}

Our objective will be to use \ref{l:extension} 
to identify the group extension
\[
1\to\Phi ({\mathcal G}_F)\to {\mathcal G}_F\to E\to 1
\]
in terms of a subspace $I_F\subset H^2(E)$. To do this we require
some basic notions from \textit{K}-theory (see \cite{M}). 

Let 
\[
l\colon \dot F/\dot F^2\to k_1F
\]
denote the canonical isomorphism between $\dot F/\dot F^2$ written
multiplicatively and additively (which is $k_1F$ by definition).
Let ${\mathbb F}_2 [k_1F]\cong H^*(E)$
denote the polynomial algebra
generated by the vector space $k_1F$, where the generators
are assumed to be one-dimensional. Then, if $a\in \dot F/\dot F^2$,
$a\ne 1$, the element $l(a)l(1-a)$ can be thought of as a quadratic polynomial
in this ring. We define $I_F\subset H^2(E)$ as the subspace generated
by these classes. If $\{e_i \mid i\in I\}$ is a basis for $\dot F/\dot F^2$,
then by definition we have that Milnor \textit{K}-theory mod 2 is given by
\[
k_*F= {\mathbb F}_2 [e_i \mid i\in I]/I_F.
\]
We will now show that our description of ${\mathcal G}_F$
as a central extension of elementary abelian groups can be
made explicit in terms of $I_F$. We use the notation
from \ref{l:extension}. 

\begin{theorem}
\label{t:central}
For the central extension 
$$1\to \Phi({\mathcal G}_F)\to {\mathcal G}_F\to E\to 1$$
$\delta (H^1(\Phi ({\mathcal G}_F)))=I_F$
and hence the \textit{W}-group
${\mathcal G}_F$ is uniquely determined by the subspace $I_F\subset H^2(E)$. In
particular there exists a basis $\{e_i\mid i\in I\}$ for $H^1(\Phi
({\mathcal G}_F))$ such that all the $k$-invariants $\delta (e_i)$ are of
the form $u_iv_i$ for $u_i,v_i\in H^1(E)$. 
\end{theorem}

\begin{proof} We make use of the diagram of extensions
described in section~\ref{s:w-groups}. First we recall a theorem due to
Merkurjev \cite{Mer}, namely that the map $H^2(E)\to H^2(G_F)$ is surjective, 
where
the kernel can be described as the subspace $I_F$ generated by the
defining relations in Milnor \textit{K}-theory. 
We claim that there is also an isomorphism 
$H^1(\Phi ({\mathcal G}_F))\cong H^1(G_F^2)^E$. Indeed, let 
$u\colon G_F^2\to {\mathbb Z}/2$ and assume it is $E$ and hence
$G_F$-invariant. 
This means that for all $g\in G_F$, $r\in G_F^2$, we have
$u(g^{-1}r^{-1}gr)=1$, hence $u(G_F^4[G_F^2,G_F])=1$ and so it defines
a unique element in $H^1(\Phi ({\mathcal G}_F))=H^1(G_F^2/G_F^4[G_F^2,G_F])$.
 From the above we obtain a diagram of exact sequences, where
the kernels and middle terms are mapped isomorphically:
\[
\xymatrix{0 \ar[r] & H^1(\Phi ({\mathcal G}_F)) \ar[r]^{\delta} \ar[d] &  H^2(E)
\ar[r] \ar[d] & H^2({\mathcal G}_F)\cap {\mathcal R} \ar[d] \ar[r] & 0 \\
0 \ar[r] & H^1(G_F^2)^E \ar[r]^{\delta'} & H^2(E) \ar[r] & H^2(G_F)
\ar[r] & 0. \\
}
\]
Hence we conclude that by choosing the right basis, we can
assume that the $k$-invariants for ${\mathcal G}_F$ are precisely the
generators for the relations in Milnor \textit{K}-theory. However
we have already mentioned that the relations there are of the form
$l(a)l(1-a)$ for $a\in\dot F/\dot F^2$, $a\ne 1$,
where $l\colon \dot F/\dot F^2 
\to k_1F$ is the canonical isomorphism. Making the
identification with the quotient of the cohomology of $H^*(E)$
we see that the $k$-invariants must indeed be products
of linear forms, completing the proof. 
\end{proof}

\begin{remark}

Actually this result only requires the injectivity part of Merkurjev's
theorem (the proof of the Milnor conjecture in dimension 2). 
Observe that for any $W$-group ${\mathcal G}_F$ we obtain all
non-trivial
$k$-invariants of the form $l(a)l(1-a)$, $a\ne 1$, because they
correspond to ${\mathbb Z}/4$ and $D_4$ quotients of ${\mathcal G}_F$.
Hence we see that all K-theoretic relations are present in
$im~\delta$, and from the injectivity of Merkurjev's theorem
and the diagram used before, we see that these provide all the
defining $k$-invariants for the $W$-group.
\end{remark}

Note that both ${\mathcal G}_F$ and $k_*F$ are unambiguously determined
by the same ideal $I_F$. From this and the main result in
\cite{MS1} we obtain

\begin{corollary}
\label{t:corollary}
${\mathcal G}_F\cong {\mathcal G}_L$ if and only if
$k_*F\cong k_*L$. Moreover if $WF\cong WL$ then
$k_*F\cong k_*L$, and the converse holds provided
$F$ and $L$ have the same level if $\langle 1,1\rangle_F$ is universal.
\end{corollary}

The following basic fact follows
from our description of the $k$-invariants
(the key idea is due to W. Gao and is used in his thesis \cite{G}, pg. 51):   

\begin{theorem}
\label{t:cohen-macaulay}
Let $F$ be a field which is not formally real
and suppose that $|\dot F/\dot F^2|=2^n$.
Then there exists a collection of 2-dimensional classes
$\zeta_1,\dots ,\zeta_r$ in $H^*({\mathcal G}_F)$ where 
$r=n + {\binom{n}{2}}- \dim H^2(G_F)$, such that
$H^*({\mathcal G}_F)$ is free and finitely generated
as a module over the subalgebra 
${\mathbb F}_2 [\zeta_1 ,\dots ,\zeta_r]$.
\end{theorem}

\begin{proof}
We consider the Lyndon-Hochschild-Serre spectral sequence
for the extension
$$1\to \Phi ({\mathcal G}_F)\to {\mathcal G}_F\to E\to 1.$$
If $H^*(E)\cong {\mathbb F}_2 [x_1,\dots ,x_n]$ and
$H^*(\Phi ({\mathcal G}_F))\cong {\mathbb F}_2 [e_1, \dots ,e_r]$
then we have already remarked that we can assume
$d_2(e_i)=u_iv_i$ where the $u_i,v_i$ are 1-dimensional linear
forms. We claim that the squares
$e_i^2$ are permanent cocycles in the spectral sequence.
For this it suffices to note that transgressions commute
with Steenrod squares, and hence
\[
d_3(e_i^2)=d_3(Sq^1(e_i))=Sq^1d_2(e_i)=
Sq^1(u_iv_i)=u_i^2v_i + u_iv_i^2=0
\]
in $E_3^{3,0}$, as we obtain an element in the ideal
generated by the image of $d_2$. 

Hence the elements $e_i^2$ are in the image of the
edge homomorphism, which in this case simply means
that they are in the image of $\Res^G_{\Phi}$. We can
therefore find polynomial classes $\zeta_i\in H^2(G)$, 
$i=1,\dots ,r$,
such that $\Res^G_{\Phi}(\zeta_i)=e_i^2$ and invoking
lemma~\ref{l:varieties} we conclude that $H^*(G)$ is free and finitely
generated over the polynomial subring which they generate.
\end{proof}

In \cite{BC} it was shown that the Poincar\'e series for
a Cohen-Macaulay cohomology ring has a very special
form. Namely, if $H^*(G)$ is free and finitely generated
over a polynomial subring ${\mathbb F}_2 [u_1,\dots ,u_k]$ 
where $u_i\in H^{n_i}(G)$, then it is of the form
$q(t)/(1-t^{n_1})\dots (1-t^{n_k})$, where $q(t)$ is a
``palindromic polynomial''
(i.e.\ if $q(t)=1+a_1t+\dots +a_dt^d$, then
$a_i=a_{d-i}$),
with integral coefficients and of degree 
$n_1+\dots + n_k - k$. Applying this we obtain
\begin{corollary}
\label{c:q_F}
If $F$ is not formally real, 
the Poincar\'e series for $H^*({\mathcal G}_F)$ is of the
form
\[
p_F(t)=\frac{q_F(t)}{(1-t^2)^r}
\]
where $q_F(t)$ is a palindromic polynomial of degree
$r=n + {\binom{n}{2}}-\dim H^2(G_F)$ 
in ${\mathbb Z} [t]$.
\end{corollary}

The proof above can be combined with \ref{t:non2C} to
establish a modified version of the above for formally real fields. 
The appeal to lemma~\ref{l:varieties} can be replaced by an appeal to
lemma~\ref{l:regular-sequence}.

\begin{theorem}
\label{t:not-quite-CM}
If $F$ is a formally real field with $|\dot F/\dot F^2|=2^n$,
then there exists a regular sequence of $2$-dimensional elements
in $H^*({\mathcal G}_F)$ 
of length $r=n + \binom{n}{2} - \dim H^2(\Gal(F_q/F))$,
i.e.\ $H^*({\mathcal G}_F)$ has depth at least equal to $r$;
this is one less than its Krull dimension.
\end{theorem}

The point here is that we have a sequence of elements restricting
non-trivially to the maximal 2-torus in the center. Hence we
conclude that $H^*({\mathcal G}_F)$ has depth at least 
one less than the Krull dimension. We cannot conclude
that it is Cohen-Macaulay, although that may occur, as we shall
see in \S 5.

\begin{remark} 
\label{r:PS-well-defined}
There is a factorization for the Poincar\'e series in the formally
real case; one factor is the expected $(1-t^2)^{-r}$, but the other
factor need not be a polynomial.  The existence of this
factorization follows from the fact \cite[10.3.4]{E} that if $\zeta_1,
\ldots, \zeta_r$ 
is a regular sequence in $H^*({\mathcal G}_F)$, then $H^*({\mathcal
G}_F)$ is a free module over ${\mathbb F}_2[\zeta_1, \ldots,
\zeta_r]$.  The
other factor to which we have referred will be given a very
concrete form later, as it is the equivariant cohomology of our
``topological model'' (see theorem~\ref{t:X_F}). 
\end{remark}

We will now show that the mod 2 cohomology of a \textit{W}-group contains
at least as much information as the Galois Cohomology
of $F$.
Let $\bar{F}$ denote a separable closure of $F$, and
let ${\mathbb G} = \Gal (\bar{F}/F)$ denote the Galois group
of this extension. Then by definition, the \textsl{Galois cohomology}
of $F$, denoted $H^*(F)$, is the cohomology of the group
${\mathbb G}$, where, as remarked earlier, we will assume mod 2
coefficients.  
Let $E={\mathbb G}/\Phi_2 ({\mathbb G})$ where $\Phi_2({\mathbb G})$
is the 2-Frattini subgroup of ${\mathbb G}$.
Note that
$H^1(F)\cong \Hom (E,{\mathbb Z}/2)$.
The quotient map above factors through
the \textit{W}-group associated to $F$ and so we have a commutative
diagram
\[
\xymatrix{{\mathbb G} \ar[r] \ar[rd] & {\mathcal G}_F \ar[d] \\
&E \\
}
\]
where by construction the group $E$ can also be identified
with the quotient
${\mathcal G}_F/\Phi ({\mathcal G}_F)$. 
Hence we have isomorphisms
$H^1(F)\cong H^1({\mathcal G}_F)\cong H^1(E)$.
Let ${\mathcal R}\subset H^*({\mathcal G}_F)$ denote the subring generated
by one dimensional classes. We have the following fundamental result

\begin{theorem}
\label{2.3}
The map $\phi \colon {\mathbb G}\to {\mathcal G}_F$ induces an isomorphism
${\mathcal R}\cong H^*(F)$.
\end{theorem}
\begin{proof}
The group $E$ is an elementary abelian
group such that $H^1(E)\cong \dot F/\dot F^2$. In cohomology
the surjection ${\mathbb G}\to E$ gives rise to a ring map
$H^*(E)\to H^*(F)$. On the other hand we have seen that
Milnor
\textit{K}-theory  
$k_*F$ can be identified with the quotient of the polynomial
ring $H^*(E)={\mathbb F}_2 [k_1F]$ by the ideal $I_F$ generated
by 2-dimensional relations,
of the form $l(a)l(1-a)$, where $l\colon \dot F/\dot F^2\to k_1F$
is the canonical isomorphism (see \cite[p.~319]{M}). In the finite case,
the \textit{K}-theory is simply expressible as a 
ring ${\mathbb F}_2 [x_1,\dots ,x_n]/(\mu_1,\dots ,\mu_s)$, where
$\mu_1,\dots ,\mu_s$ are an irredundant collection of quadratic 
polynomials in the variables $x_1,\dots ,x_n$. 
By Voevodsky's theorem (the Milnor conjecture \cite{V}), we can
identify $k_*F$ with $H^*({\mathbb G})$ using the map
$H^*(E)\to H^*({\mathbb G})$ i.e.\ we can identify its kernel
with the ideal $I_F$.

We want to fit the \textit{W}-group into this picture. Note that as a consequence
of our previous remarks, ${\mathcal R}$ surjects onto $H^*(F)$. Now we have
already established in \ref{t:central} that
if the group ${\mathcal G}_F$
is expressed as a central extension

$$1\to \Phi ({\mathcal G}_F)\to {\mathcal G}_F\to E\to 1$$
defined by $k$-invariants 
$\rho_j\in H^2(E)$, $j\in \Omega$, then we can assume
that the ideal $J$ generated
by the $\rho_j$ is precisely the ideal $I_F$
generated by the basic Steinberg relations
in \textit{K}-theory. 
Hence we can express ${\mathcal R}$ as a quotient of 
$H^*(E)/I_F$.
Putting this together we deduce that we have a factorization
of the isomorphism in the Milnor conjecture into a composition
of two epimorphisms:
\[
k_*F\longrightarrow {\mathcal R}\longrightarrow H^*(F)
\]
from which the result follows.
\end{proof}

One can think of the \textit{W}-group ${\mathcal G}_F$ as the 
central extension of an elementary abelian group by an elementary
abelian group such that its associated $k$-invariants are precisely
the generating
relations for \textit{K}-theory.
Note that it is a consequence of
the Milnor Conjecture that ${\mathcal R}$ is precisely the quotient
$H^*(E)/I_F$. Computing
$H^*({\mathcal G}_F)$ using the Lyndon-Hochschild-Serre
spectral sequence associated to the
extension above, this means that 
$$E_{\infty}^{*,0}=E_3^{*,0}
\cong H^*(F)\cong {\mathcal R}$$
i.e.\ no differentials beyond $d_2$ can hit the cohomology of the
base. We will see later that a general
collapse of these spectral sequences
at the $E_3$ term for all \textit{W}-groups is not altogether unlikely.

We can use the above to deduce the following 

\begin{corollary}
\label{2.4}
If $F$ is a field of characteristic different from 2
such that ${\mathcal G}_F\ne{\mathbb Z}/2$, then the following four
conditions are equivalent:
\begin{enumerate}

\item every $a\in H^1({\mathcal G}_F)$ is nilpotent

\item every positive dimensional element in $H^*(F)$ is
nilpotent,

\item $F$ is not formally real,

\item ${\mathcal G}_F$ satisfies the 2C condition.
\end{enumerate}
\end{corollary}

\begin{proof}
The equivalence of (1) and (2) follows from theorem~\ref{2.3}. The
equivalence of 
(2) and (3) is theorem 1.4 in \cite{M}. The equivalence of (3) and (4)
is proposition~\ref{t:2C}.
\end{proof}

In \cite{AK} it was shown that if $G$ is a finite 2-group satisfying the
2C condition, then it contains non-zero cohomology classes
$x\in H^*(G,{\mathbb F}_2)$ such that $\Res^G_H(x)=0$ for every proper
subgroup $H\subset G$ (these are called
\textsl{essential} cohomology classes). We conclude 

\begin{proposition}
\label{p:essential-exterior}
If $F$ is not formally real and $|\dot F/\dot F^2|<\infty$,
and $|{\mathcal G}_F|>2$,
then $H^*({\mathcal G}_F)$ contains essential cohomology classes.
Furthermore if $x$ is an essential class then it is exterior, i.e.
$x^2=0$.
\end{proposition}

\begin{proof}
Only the second part needs justification.
Given our hypotheses we have a non-trivial relation
in \textit{K}-theory, i.e.\ $I_F\ne 0$.
Hence we see that there 
exist non-zero classes
$x_1,x_2\in H^1({\mathcal G}_F)$ 
such that $x_1x_2=0$.
On the other hand due to the fact that $x$ is essential, it must be
divisible by any one dimensional class. Choose $u_1,u_2$ such that
$x=x_1u_1=x_2u_2$; then obviously $x^2=x_1u_1x_2u_2=x_1x_2u_1u_2=0$.
\end{proof}

\begin{remark}
Observe that the same proof shows that all products of essential
classes are zero. In other words, the ideal of essential classes
${\mathcal E}\subset H^*({\mathcal G}_F)$ is
a ring with trivial products. An interesting problem
would be to obtain a description of ${\mathcal E}$. For example, one can
ask for the minimal degree of a homogeneous element $x\in {\mathcal E}$.
As an example consider the case when $F$ is a local field 
(a finite extension of some
${\mathbb Q}_p$) then we know (see \cite{Se}, 4.5) that the cup product
$$H^1({\mathcal G}_F)\times H^1({\mathcal G}_F)\to {\mathbb Z}/2\subset 
H^2({\mathcal G}_F)$$
is a non--degenerate  bilinear pairing (here we use the fact that
the Galois cohomology is a 2-dimensional Poincar\'e duality
algebra). Hence there is a fixed 
element
$a\in H^2({\mathcal G}_F)$ such that if $x,y\in H^1({\mathcal G}_F)$,
then $xy=0$ or $xy=a$. Clearly $a\in {\mathcal E}$, and so the essential ideal
has elements of degree two. 
\end{remark}

It remains to consider the case when $F$ is formally real.
As we have seen, the subgroup
$\Phi ({\mathcal G}_F)\subset {\mathcal G}_F$
is a central elementary abelian subgroup of index 2 in any maximal
elementary abelian subgroup in ${\mathcal G}_F$. Consequently there exists
a non-nilpotent class $\alpha\in H^1({\mathcal G}_F)$. We can make explicit
the identification between $H^1({\mathcal G}_F)$ and $\dot F/\dot F^2$
as follows. If $a\in \dot F$, then we define a homomorphism
${\mathcal G}_F\to {\mathbb Z}/2$ by the formula
$$x_a(\sigma )=(\sqrt{a})^{\sigma}/\sqrt{a}, ~~~\sigma\in {\mathcal G}_F.$$
In the sequel this class will often be denoted
by $[a]$. In terms of this we can now specify a non-nilpotent
class when the field is formally real. 

\begin{theorem}
Let $F$ denote a field of characteristic different from two and
such that $|\dot F/\dot F^2|<\infty$.
Then $F$ is formally real if and only if $x_{-1}\in H^1({\mathcal G}_F)$
is non-nilpotent.
\end{theorem}

\begin{proof}
Assume that $F$ is formally real. As we saw in the proof of \ref{t:2C},
we can choose an involution $\sigma\in {\mathcal G}_F - \Phi ({\mathcal G}_F)$
such that $(\sqrt{-1})^{\sigma}=-\sqrt{-1}$; hence we deduce that 
$x_{-1}(\sigma)\ne 1$. This means of course that $x_{-1}$ restricts
non trivially to $H^1(\langle\sigma\rangle)$, and hence $x_{-1}$ is non-nilpotent.
The converse has already been established.
\end{proof}

\begin{remark}
Actually the condition $|\dot F/\dot F^2|<\infty$ can be removed 
if one is willing to apply
the other half of the Milnor Conjecture, namely the isomorphism
$$\phi_n:H^n(F)\to I^n/I^{n+1}$$
given via 
$$\phi_n([a_1]\cup \dots\cup [a_n])=
[\langle\langle -a_1,\dots ,-a_n\rangle\rangle]$$
where $I^n$ is the n-th power of the fundamental ideal in $WF$
and 
$$[\langle\langle -a_1,\dots ,-a_n\rangle\rangle]\in I^n/I^{n+1}$$
is the class
of the Pfister form (see \cite{V}, \cite{OVV} and \cite{Lam1}).
In particular we have that $\phi_n ([-1]^n)=
[\langle\langle 1,\dots ,1\rangle\rangle]$; however from the 
Arason-Pfister Theorem (see \cite{Lam1}, Ch.10, Section 3, Th. 3.1)
we see that $[-1]^n=0$ in $H^*(F)$ if and only if 
$\langle\langle 1,\dots ,1\rangle\rangle$ is isotropic. Since
this form is \emph{never} isotropic over a formally real field, we
see that $[-1]$ is non--nilpotent in $H^*(F)$ and consequently
in $H^*({\mathcal G}_F)$.
\end{remark}

In the case of a field $F$ which is not formally real, we 
can specify exactly the height of the nilpotent class $[-1]$.

\begin{proposition}
Let $F$ be a field that is not formally real, and let $n$ denote
the smallest positive integer such that $[-1]^n=0$ in $H^*({\mathcal G}_F)$.
Then, if $s(F)$ denotes the level of the field $F$, $s(F)=2^{n-1}$. 
\end{proposition}
\begin{proof}
Clearly the height of $[-1]$ in $H^*(F)$ is the same as its height
in $H^*({\mathcal G}_F)$. However from the previous remark we see
that $[-1]^n=0$ if and only if $\langle\langle 1,\dots , 1\rangle\rangle$
is isotropic, and therefore hyperbolic (see \cite{Lam1}, Chapter 10).
Therefore the Pfister form $\langle\langle 1,\dots ,1\rangle\rangle$
($n-1$ times) represents $-1$, and this means that $-1$ is expressible
as a sum of $2^{n-1}$ squares, whence our proposition follows.
\end{proof}

This type of argument can also be used to find an upper bound on
the highest degree in which ${\mathcal R}\subset H^*({\mathcal G}_F)$
has non--zero classes (we denote this by $c(F)$). 

\begin{theorem}
Let $F$ be a field that is not formally real and such that
$|\dot F/\dot F^2|=2^n$; then $c(F)\le n$.
\end{theorem}
\begin{proof}
 From Kneser's Theorem (see \cite{Lam1}, Chapter 11, Th.4.4), we see
that each Pfister form 
$\langle\langle -a_1,\dots ,-a_{n+1}\rangle\rangle$
with $a_1,\dots ,a_{n+1}\in \dot F$ is hyperbolic; hence each
product $[a_1]\cup\dots \cup [a_{n+1}]\in H^{n+1}({\mathcal G}_F)$
is zero.
\end{proof}
\end{section}

\begin{section}{Topological Models}
\label{s:top-models}

As a consequence of the cohomological analysis carried out in the previous
section, it is fairly evident that the key step in computing the cohomology
of a \textit{W}-group is calculating the quotient algebra 
$H^*({\mathcal G}_F)/(\zeta_1,\dots ,\zeta_r)$.
In this section we will construct a topological space $X_F$ together with
a very explicit action of $E_n=({\mathbb Z}/2)^n$ on it, such that the quotient
above can be identified with the equivariant cohomology $H^*_{E_n}(X_F)$.
This \textsl{topological model} has geometric properties which reflect
the field theory in a very nice way, and as we shall see, computing
the equivariant cohomology rings is an interesting problem in its own
right. 

Before proceeding with the description of our construction and
proving its main properties, we briefly recall some
basic facts about equivariant cohomology that will be used
in the rest of the paper.  All of these facts can be found
in standard texts, e.g.\ \cite[Ch.~V]{AM}, \cite{AP},
\cite{Bredon}.  

Let $X$ be a topological space with the action of a finite group $G$ on
it. Let $EG$ denote a universal $G$--space i.e. a contractible
space with a free action of $G$.
The {\sl Borel construction} is defined as the orbit space
$EG\times_G X = (EG\times X )/G$, where $G$ acts diagonally on
$EG\times X$. The cohomology of this construction is
denoted by $H^*_G(X) = 
H^*(EG \times_G X,{\mathbb F}_2)$, and is referred to  
as the mod 2 \textsl{equivariant
cohomology} of $X$. For our applications we will always assume that
$X$ is a compact space, in most cases a closed manifold, and that
$G$ is a finite $2$-group. As the
action on $EG$ is free, the projection
map $EG\times X\to EG$ gives rise to a fibration $EG\times_G X\to BG$,
where $BG=EG/G$ is the classifying space of $G$, and the fiber is
$X$. Similarly if the $G$-action on $X$ is free, we get a fibration
$EG\times_G X\to X/G$ with contractible fiber $EG$, hence we have
a homotopy equivalence $EG\times_G X\simeq X/G$.

The basic structural result concerning
equivariant cohomology is that 
$H^*_G(X)$ has Krull dimension equal
to the largest rank of a $2$-elementary abelian isotropy subgroup.
We will use two simple consequences of this:

\begin{lemma}
If $X$ is a finite dimensional $G$-space,
then $H^*_G(X)$ is finite dimensional if and only if $G$
acts freely on $X$, in which case $H^*_G(X)\cong H^*(X/G)$. 
\label{ec:finite-dimension}
\end{lemma}

\begin{lemma}
If the isotropy groups $G_x$ are cyclic (or trivial) for all $x
\in X$, then $H^*_G(X)$ is eventually periodic. \label{ec:periodic}
\end{lemma}

The \textsl{singular set} $\Sing_G(X)$  is the set of
points of $X$ which have non-trivial isotropy groups, or equivalently,
the set of points of $X$ which are not permuted freely by $G$.
In sufficiently large dimensions $H^*_G(X)$ is determined by the
equivariant cohomology of the singular set. More precisely, we have
\begin{lemma}
If $i > \dim X$, then the inclusion $\Sing_G(X) \subset X$ induces an 
isomorphism $H^i_G(X) \cong H^i_G(\Sing_G(X))$.
\label{ec:large-dimensions}
\end{lemma}

We will also use the following simple fact:

\begin{lemma}
If $H \subset G$ is a subgroup, $H^*_G(G/H \times Y) \cong H^*_H(Y)$.
\label{ec:frobenius} 
\end{lemma}

The inclusion $X \subset EG \times_G X$ gives a map $H^*_G(X) \to
H^*(X)$, and if this map is surjective, we say that
the cohomology of $X$ is \textsl{totally nonhomologous to zero} in the
cohomology of the \textsl{Borel construction} $EG \times_G X$. 
Writing the Poincar\'e series for a graded ring $R$ as $P(R)$ we have:

\begin{lemma}
\label{ec:TNHZ}
If $H^*(X)$ is totally nonhomologous to zero in $H^*_G(X)$, then 
$P(H^*_G(X))
= P(H^*(X)) \cdot P(H^*(G))$.
\end{lemma}

We now state and prove the main result in this section.

\begin{theorem}
\label{t:X_F}
Let $F$ denote a field of characteristic different from 2
such that $|\dot F/\dot F^2|=2^n$. 
Let $r$ denote the rank of $\Phi ({\mathcal G}_F)$. There exists
a homomorphism $\rho\colon ({\mathbb Z}/2)^n\to O(2)^r$ which defines an action
of $E_n=({\mathbb Z}/2)^n$ on $X_F\cong (\Sone)^r$ with the 
following properties:
\begin{enumerate}
\item $E_n$ acts freely on $X_F$
if and only if $F$ is not formally
real.
\item $E_n$ acts with non-trivial cyclic isotropy subgroups if and
only if $F$ is formally real.
\item There exists a regular sequence $\{\mu_F\}\subset 
H^2({\mathcal G}_F)$ of length 
$$r=n+{\binom{n}{2}} - \dim 
H^2(\Gal(F_q/F))$$ such that
$$H^*({\mathcal G}_F)/(\mu_F)\cong H^*_{E_n}(X_F).$$
\end{enumerate}
\end{theorem}

\begin{proof}
To prove this we will make use of the
fact that ${\mathcal G}_F$ is defined as a central extension
$$1\to \Phi ({\mathcal G}_F)\to {\mathcal G}_F\to E_n\to 1$$
where we can assume that the $r$ $k$-invariants are of the
form 
$$u_1v_1,\dots ,u_rv_r,$$ 
for elements $u_s,v_t\in H^1(E_n)$.
Given any element $\beta\in H^1(E_n)$, there is a representation
$\rho_{\beta}\colon E_n\to O(1)$ such that it has
Stiefel-Whitney class $w_1(\rho_{\beta})=\beta$.
Hence we may construct a direct sum representation
$\rho_{u_t}\oplus \rho_{v_t}=\rho_{u_tv_t}\colon E_n\to O(2)$ such that 
it has second Stiefel-Whitney class
$w_2(\rho_{u_tv_t})=w_1(\rho_{u_t})w_1(\rho_{v_t})=u_tv_t$. Taking
a product we obtain a homomorphism $\rho \colon E_n\to O(2)^r$ such that
on each factor the corresponding $w_2$ is mapped to the appropriate
$k$-invariant.

Consider the action of $O(2)$ on $\Sone$ induced by matrix 
multiplication on ${\mathbb R}^2$. Restricted to the subgroup
$O(1)\times O(1)\subset O(2)$ we obtain the following action:
let $x_1,x_2$ denote generators of $G=O(1)\times O(1)$, then if
$z\in\Sone$, $x_1(z)=\bar{z}$, $x_2(z)=-\bar{z}$. Here
we have used complex conjugation to simplify the notation. 
Note that the subgroups $\langle x_1\rangle$, $\langle x_2\rangle$ 
are isotropy subgroups,
while $x_1x_2$ acts freely via multiplication by $-1$.
We can use
the universal $G$ space $EG$ to obtain
the usual fibration
\[
\xymatrix{\Sone \ar[r] & EG\times_G \Sone \ar[d] \\
& BG
}
\]
If $e\in H^1(\Sone)$ is a generator and 
$H^*(G)\cong{\mathbb F}_2 [x_1, x_2]$ (a slight abuse of notation),
then the transgression can be computed as
$d_2(e)=x_1x_2$. This follows from the fact that 
$d_2(e)$ restricts trivially to the cohomology
of both of the isotropy subgroups
(indeed both restricted fibrations have sections) and non-trivially
to the cohomology of $\langle x_1x_2\rangle$. 
In fact $\pi_1(EG\times_G\Sone)
\cong {\mathbb Z}/2*{\mathbb Z}/2$, the infinite dihedral group.
Hence in the fibration 
\[
\xymatrix{\Sone \ar[r] & EO(2)\times_{O(2)} \Sone \ar[d] \\
&BO(2)
}
\]
the cohomology generator on the fiber will transgress to the
Stiefel-Whitney class $w_2=x_1x_2$. 

Now let $X_F$ denote the space $(\Sone)^r$ with the action
of $E_n$ through the homomorphism 
$\rho \colon E_n\to O(2)^r$ described previously.
Using the naturality of the 
corresponding fibration for the Borel construction
$X_F\times_{E_n}EE_n$, one can verify that the cohomology
generators
on the fiber will transgress to the images of the $w_2$ classes
under the map induced by $\rho$, which by construction are
precisely the desired $k$-invariants $u_1v_1,\dots, u_rv_r$.
Note that $X_F$ is
a product of $r$ $E_n$-spaces, as the action on each
coordinate is independent of the others. Let
$Y_F = X_F\times_{E_n}EE_n$, then we have the usual fibration
\[
\xymatrix{X_F \ar[r] & Y_F \ar[d] \\
	   &BE_n
}
\]
which corresponds to a group extension

$$1\to L\to \Upsilon_F\to ({\mathbb Z}/2)^n\to 1$$
where $\Upsilon_F = \pi_1 (Y_F)$ and $L$ is a free abelian group
of rank $r$. Note that if the $E_n$ action is free,
then $Y_F$ is an
aspherical, compact $r$-dimensional manifold, with necessarily
torsion-free fundamental group. 

We now consider the mod 2
spectral sequence for this extension (or equivalently the
fibration). Let $H^*(E_n)\cong {\mathbb F}_2[x_1,\dots ,x_n]$ and
$H^*(X_F)=\Lambda (e_1, \dots ,e_r)$. We have seen
that \emph{by construction} the fibration
will have $k$-invariants $d_2(e_i)=u_iv_i$, for $i=1,\dots ,r$.

The subgroup $2L\subset\Upsilon_F$ is a normal subgroup, and
hence fits into a commutative diagram of extensions:
\[
\xymatrix{
 &  1 \ar[d] & 1 \ar[d] \\
& 2L \ar@{=}[r] \ar[d] & 2L \ar[d] \\
1 \ar[r] & L \ar[r] \ar[d] & \Upsilon_F \ar[r] \ar[d] & ({\mathbb Z}/2)^n \ar[r] \ar@{=}[d] &	1 \\
1 \ar[r] & ({\mathbb Z}/2)^r \ar[r] \ar[d] & {\mathcal G}_F^{\prime} \ar[r]
\ar[d] & ({\mathbb Z}/2)^n \ar[r] & 1\\ 
& 1 & 1 \\
}
\]
Now observe that in the spectral sequence for the bottom row,
the transgressions are precisely $\{u_iv_i\mid i=1,\dots ,r\}$.
As the
original action on $L$ was a sum of rank one sign-twists, we
conclude that the quotient group ${\mathcal G}_F^{\prime}$ is expressed
as a central extension of $({\mathbb Z}/2)^n$ by $({\mathbb Z}/2)^r$ with
precisely the same extension class as ${\mathcal G}_F$. We conclude
from \ref{l:extension} that ${\mathcal G}_F^{\prime}\cong {\mathcal G}_F$. 
We now make use of
the 2 columns; comparing the spectral sequences for them
we see that we may choose $\tilde{e}_i$ generating $H^1(2L)$
such that the $d_2(\tilde{e}_i)$ restrict to the squares of the
polynomial generators in $H^*(({\mathbb Z}/2)^r)$.
These elements, which we denote by $\zeta_i$, $i=1,\dots ,r$,
form a regular sequence in $H^*({\mathcal G}_F)$ by
lemma~\ref{l:regular-sequence}.   
As a consequence of this, the spectral sequence collapses
at $E_3$ and 
\[
H^*(Y_F)\cong H^*(\Upsilon_F)\cong H^*({\mathcal G}_F)/(\zeta_i).
\]

We now consider the case when $F$ is not formally real. As we
have seen, this is equivalent to having ${\mathcal G}_F$ satisfy
the 2C condition. In other words, $\Phi ({\mathcal G}_F)$ is  
the unique elementary abelian subgroup of maximal rank, and
it is central. This means that no cyclic subgroup in the quotient
$E_n$ can split off; cohomologically this means that if we
consider the restricted extension
$$1\to \Phi ({\mathcal G}_F)\to H_C\to C\to 1$$
where $C\subset E_n$ is any cyclic subgroup, then the restricted
$k$-invariants must generate $H^2(C)={\mathbb F}_2$. Looking at the 
associated fibration for the $C$-action on $X_F$; we see that
the $k$-invariants must also necessarily generate $H^2(C)$; hence
$C$ acts freely on $X_F$ and we conclude that the entire group
$E_n$ must be acting freely.

It remains to see what happens in the formally real case. Here we
know that $\Phi ({\mathcal G}_F)$ is a central elementary abelian subgroup,
of index 2 in any maximal elementary abelian subgroup. Moreover,
if $s_1, s_2$ are commuting
involutions in ${\mathcal G}_F-\Phi ({\mathcal G}_F)$, then
$s_1s_2\in \Phi ({\mathcal G}_F)$. If there were a subgroup $A\cong 
{\mathbb Z}/2\times {\mathbb Z}/2$ in $E_n$ which fixed a point in $X_F$, then
the associated splitting
would imply the existence of an elementary abelian subgroup
of rank equal to two more than the rank of $\Phi ({\mathcal G}_F)$,
a contradiction. On the other hand we do know that any cyclic
summand in $E_n$ which forms a maximal elementary abelian subgroup
together with $\Phi ({\mathcal G}_F)$ must split, hence it splits
cohomologically and must have a fixed point.
\end{proof}

If $F$ is not formally real, we have shown that $E_n$ acts freely
on $X_F$. Hence we have a homotopy equivalence 
$$Y_F\simeq X_F/E_n$$
and so $H^*_{E_n}(X_F)\cong H^*(X_F/E_n)$. Now $X_F/E_n$
is a compact $r$-dimensional manifold, such that its mod 2 Poincar\'e
series is precisely the polynomial $q_F(t)$ defined in \ref{c:q_F}.
In the formally real case the the equivariant cohomology
is infinite dimensional, although eventually periodic, as the
isotropy is cyclic (\ref{ec:periodic}). An interesting alternative here would be
to compute the equivariant Tate cohomology  
$\widehat{H}_{E_n}^*(X_F)$; this
invariant vanishes for free actions and more generally
agrees with the ordinary equivariant cohomology in sufficiently
high dimensions. For the case of fields which are formally real,
the entire Tate Cohomology is periodic  
and it would
seem possible to find an expression for it in field-theoretic
terms, but we will not explore this any further here. 

The result above shows that the calculation of $H^*({\mathcal G}_F)$ can
be reduced to computing the cohomology of a crystallographic group.
When the field is not formally real then in fact it reduces to 
computing the cohomology of a Bieberbach group. It is worthwhile
to compare the spectral sequences arising from the two extensions
\[
1\to M_F\to\Upsilon_F\to E_n\to 1
\]
\[
1\to\Phi({\mathcal G}_F)\to {\mathcal G}_F\to E_n\to 1
\]
using the
morphism of extensions which we defined above. The $E_2$ terms
are mapped via
\[
{\mathbb F}_2[x_1,\dots ,x_n]\otimes {\mathbb F}_2[e_1,\dots ,e_r]
\longrightarrow {\mathbb F}_2[x_1,\dots ,x_n]\otimes \Lambda
(e_1,\dots ,e_r)
\]
where $e_i^2$ is sent to zero. Using this and the fact that the
$e_i^2$ are permanent cocycles (due to the form of the $k$-invariants),
we conclude that at any stage the map above induces an isomorphism
of spectral sequences
\[
E_r^{*,*}({\mathcal G}_F)/(e_1^2,\dots ,e_r^2)
\cong E_r^{*,*}(\Upsilon_F).
\]

We can deduce a few facts from this. First of all, the spectral sequences
will both collapse at the same stage, in particular $E_3({\mathcal G}_F)
=E_{\infty}({\mathcal G}_F)$ if and only if $E_3(\Upsilon_F)=
E_{\infty}(\Upsilon_F)$. Finding conditions that imply the collapse at $E_3$ 
for either type of extension is a basic question which also relates
to the collapse of an associated Eilenberg-Moore spectral sequence,
as we shall discuss later on.
 From the above we also can conclude that 
$E_3^{*,0}({\mathcal G}_F)=E_{\infty}^{*,0}({\mathcal G}_F)
\cong E_{\infty}^{*,0}(\Upsilon_F)$, 
hence we see that 
\[
H^*(F)\cong {\mathcal R}\subset H^*_{E_n}(X_F)
\]
where we can now identify ${\mathcal R}$ with the subring
in $H^*_{E_n}(X_F)$ generated by one dimensional elements.
For the case of fields which are not formally real this means that the
Galois cohomology $H^*(F)$ can be computed as the subring generated 
by one dimensional elements in the cohomology of a compact
manifold.

A consequence of the proof of \ref{t:X_F} is a description of the isotropy
subgroups in terms of group theoretic data for ${\mathcal G}_F$. 
The elements in $E_n$ which can fix a point in $X_F$ are exactly
those which are the homomorphic image of a (noncentral) involution
under the natural projection ${\mathcal G}_F\to E_n$. Moreover, given
any two elements which fix points in $X_F$, their product will act
freely on it. Another observation for the case of fields which are
not formally real 
is that the fact that $E_n$ acts freely on $X_F$ inducing a trivial
action in mod 2 cohomology implies (see \cite{Ca}) that
$n\le r$. This can be stated in a purely algebraic way as follows:
if $F$ is a field which is not formally real, then
$\dim H^2(\Gal(F_q/F))\le {\binom{n}{2}}$. In the formally
real case (using a modified result involving the co-rank of an 
isotropy subgroup of maximal rank) the bound is ${\binom{n}{2}}+1$.

One would like to comprehend in a precise way how unique our construction
is. For this we must keep in mind the data which was used. To begin we
fixed a choice of the $k$-invariants (or, equivalently, a basis for
$H^1(\Phi({\mathcal G}_F))$) such that each of the $r$ $k$-invariants
was a product. Then, given a $k$-invariant of the form $\alpha=uv$ we
constructed a map $E_n\to O(2)$ and from there an action 
on $\Sone$. Now $E_n$ will act on $H^1(\Sone,{\mathbb Z})\cong{\mathbb Z}$; there will
be a subgroup $H_{\alpha}\subset E_n$ which acts trivially on this
module, and which has index at most 2 in $E_n$. This subgroup corresponds
to the one dimensional cohomology class $u+v$; note that $H_{\alpha}=E_n$
if and only if $u=v$.
If the action is non-trivial, this module, denoted $M_{H_{\alpha}}$, 
fits into a short exact sequence of ${\mathbb Z}E_n$-modules
$$0\to M_{H_{\alpha}}\to {\mathbb Z} [E_n/H_{\alpha}]\to{\mathbb Z}\to 0.$$

Applying group cohomology to the sequence above, it follows
that 
$$H^2(E_n,M_{H_{\alpha}})\cong H^2(H_{\alpha},{\mathbb Z}),$$ 
and
furthermore as it is 2-torsion, 
the mod 2 reduction map $M_{H_{\alpha}}\to {\mathbb F}_2$
induces a monomorphism $0\to H^2(E_n,M_{H_{\alpha}})\to H^2(E_n,{\mathbb F}_2)$. 
Our construction can be interpreted as showing the existence of an element
$\bar{\alpha}\in H^2(E_n,M_{H_{\alpha}})$ which maps to $\alpha$ under this
map, note that it is automatically unique
once the representation
$M_{H_{\alpha}}$ has been fixed
(this also works if $H_{\alpha}=E_n$). Geometrically this corresponds to
taking the extension $\pi_1(\Sone\times_{E_n}EE_n)$, with the action inducing
$M_{H_{\alpha}}$. We do this for each $k$-invariant and then take a 
direct sum, to obtain the unique class $\alpha_F\in H^2(E_n, M_F)$ which under
mod 2 reduction maps to the sum of all of the $k$-invariants (here
$M_F=\bigoplus M_{H_{\alpha_i}}$). Again, this lifting will be unique for
a prescribed module $M_F$. However, there are many possible choices
for the module. 

As an example, consider the case of the universal \textit{W}-group
$W(2)$. The $k$-invariants are $\{x_1^2,x_2^2,x_1x_2\}$ and our construction
in this case provides $X_F=(\Sone)^3$ with a ${\mathbb Z}/2\times{\mathbb Z}/2$ action
defined by $x_1(z_1,z_2,z_3)=(-z_1,z_2,\bar{z}_3)$ and
$x_2(z_1,z_2,z_3)=(z_1,-z_2,-\bar{z}_3)$. Now we can modify the
$k$-invariants without changing the quotient group to obtain
$\{x_1^2+x_1x_2,x_2^2+x_1x_2,x_1x_2\}$. In this case our construction
provides $X_F^{\prime}\cong (\Sone)^3$ with an action given by
$x_1(z_1,z_2,z_3)=(-z_1, \bar{z}_2,-\bar{z}_3)$ and
$x_2(z_1,z_2,z_3)=(\bar{z}_1,-z_2,\bar{z}_3)$. Note that 
the first Betti number for $X_F/E_2$ is 2 whereas the first betti number
for $X_F^{\prime}/E_2$ is zero. Nevertheless the two spaces have the
same mod 2 cohomology groups.

Thus, to summarize, while the \textit{W}-group ${\mathcal G}_F$ depends only
on the subspace of $H^2(E_n)$ generated by the $k$-invariants
(lemma~\ref{l:extension}), our topological model $X_F$ is constructed
from, and depends on, a particular choice of basis for that subspace.
It follows from remark~\ref{r:PS-well-defined} that 
the cohomology groups $H^*(Y_F)$ are determined by ${\mathcal
G}_F$. 
Although different models may yield different
ring structures for $H^*(Y_F)$, 
they will all contain the Milnor
\textit{K}-theory of $F$ (theorem~\ref{2.3}) and the ring structure of
$H^*({\mathcal G}_F)$ can in principle be recovered from any one
of them by solving
an extension problem. 

The case of fields which are not formally real is
particularly interesting, 
as we are simply dealing with certain quotients of $(\Sone)^r$
by the free action of an elementary abelian 2-group.
By our construction, the subgroup
$H_F=\cap_{\alpha_F}H_{\alpha}$ is acting freely and homologically trivially
on the product of circles. In fact the elements in this group are acting
on each circle either as the identity or as multiplication by $-1$. From
this it follows that $X_F/H_F\cong (\Sone)^r$. This space,
denoted $\bar{X}_F$, now has a free action of the quotient group
$E_n/H_F$ which induces a faithful representation on $H^1(\bar{X}_F,{\mathbb Z})$.
This is the standard model for a Bieberbach group (see \cite{Ch}). In the language
of extensions, what we have done is take a \emph{maximal} abelian subgroup
$\pi_1(\bar{X}_F)\subset\Upsilon_F$, with quotient the \textsl{holonomy
group} $E_n/H_F$. Going back to the examples above for $W(2)$, we see
that in the first model the subgroup generated by $x_1x_2$ acts via
diagonal multiplication by $-1$, hence trivially in homology. Dividing
out by this action we obtain $\bar{X}_F \cong (\Sone)^3$ with 
holonomy ${\mathbb Z}/2=E_2/\langle x_1x_2\rangle$. We can locate this space
in the list of all compact connected flat riemannian 3-manifolds
(which have been classified, see \cite{W}, page 122). From the fact
that $X_F/E_2$ 
is non-orientable, with $\beta_1=2$ we deduce that this must be the
manifold of type ${\mathcal B}_2$. On the other hand, the other model
does have the full $E_2$ as holonomy group and
$X_F^{\prime}/E_2$ is orientable, with $\beta_1=0$. Looking at the
classification, we conclude that this must be the celebrated
\textsl{Hantzsche-Wendt}
manifold (of type ${\mathcal G}_6$ in Wolf's notation).  

 From the discussion above we see that the classification of \textit{W}-groups
for fields which are not formally real will involve information about
the classification of compact connected, flat Riemannian manifolds
with elementary 2-abelian holonomy.
In dimension $3$ the possibilities are fairly limited: $3$ orientable
ones (${\mathcal G}_1, {\mathcal G}_2, {\mathcal G}_6$) and $4$ non-orientable
ones (${\mathcal B}_1, {\mathcal B}_2, {\mathcal B}_3, {\mathcal B}_4$).  
Of course we are most interested in a classification ``up to mod 2
cohomology''. This provides an additional simplification: there
are only 2 distinct mod 2 Poincar\'e series which occur among the $7$
manifolds above: the ones with Poincar\'e series
$1+2t+2t^2 +t^3$ (${\mathcal G}_6, {\mathcal B}_2, {\mathcal B}_4$) and those
with Poincar\'e series $1+3t+3t^2+t^3$ (${\mathcal G}_1,{\mathcal G}_2, 
{\mathcal B}_1, {\mathcal B}_3$). The first type occurs for $W(2)$ whereas
the second kind occurs for ${\mathcal G}_F=({\mathbb Z}/4)^3$. This completely
describes the possibilities for the mod 2 cohomology of \textit{W}-groups
of formally real fields $F$ such that $\dim \Phi ({\mathcal G}_F)=3$.
Of course we can subdivide this into the two cases
\begin{enumerate}
\item $|\dot F/\dot F^2|=4$, $\dim H^2(G_F)=0$ 
and $q_F(t)=1+2t+2t^2+t^3$ 

\item $|\dot F/\dot F^2|=8$, 
$\dim H^2(G_F)=3$ and $q_F(t)=1+3t+3t^2+t^3$. 
\end{enumerate}

Making use of
this type of classification in higher dimensions would be much more
complicated, but we can still make complete statements for
dimension four.
Even though there are 74 equivalence classes of compact
4-dimensional euclidean space forms (see \cite[p.~126]{W} and \cite{H}), the
possibilities for us are severely restricted by the fact that
the holonomy can be at most $({\mathbb Z}/2)^4$ and the $k$-invariants
must be linearly independent. Hence $|\dot F/\dot F^2|$ can only
have either $8$ or $16$ elements. As before we obtain only two
possible cases: 

\begin{enumerate}
\item $|\dot F/\dot F^2|=8$, 
$\dim H^2(G_F)=2$, and $X_F/E_3$ with Poincar\'e
series equal to $q_F(t)=1+3t+4t^2+3t^3+t^4$

\item $|\dot F/\dot F^2|=16$, 
$\dim H^2(G_F)=6$, and $X_F/E_4$ with Poincar\'e
series equal to $q_F(t)=1+4t+6t^2+4t^3+t^4$. 
\end{enumerate}

To make this discussion complete, we list fields
which give rise to the two possible Poincar\'e series.
Let $K$ denote a field such that ${\mathcal G}_K=W(2)$ (this is discussed
in \S\ref{s:universal}). 
If we let $F=K((t))$, then the cohomology of ${\mathcal G}_F$ will give
rise to the Poincar\'e series in (1). 

If we take 
$F={\mathbb F}_p ((t_1))((t_2))((t_3))$ where $p$ is a prime congruent
to $3$ mod $4$, then ${\mathcal G}_F$ gives rise to the Poincar\'e
series in (2). 
An analysis of the possible cohomology
\emph{rings} is of course a more delicate issue.

We will now look at further
explicit examples of interest in field theory.

\begin{example}{$F={\mathbb Q}_2$, the field of $2$-adic numbers.}
\label{Q2}
In this case, the vectors $[-1],[2], [5]$ form a basis for $\dot F/\dot F^2$. Let
us denote by $x_{-1}, x_2, x_5$ the elements in $H^1({\mathcal G}_{\mathbb Q_2})$
which correspond to this basis. Then the $k$-invariants for this \textit{W}-group
are 
$$\{x_2^2, x_5^2, x_{-1}x_2, x_{-1}x_5, x_{-1}^2 + x_2x_5\}$$ 
(see \cite{MS1}).
To apply \ref{t:X_F}, we must find equivalent $k$-invariants expressed as products
of one dimensional classes. By a simple change of basis, we can modify the list
above to get 
$$\{x_2^2,x_5^2, x_{-1}x_2, x_{-1}x_5, (x_{-1}+x_2)(x_{-1}+x_5)\}$$
without altering the group (in fact this simply corresponds to making the
right choice of generators for $H^1(\Phi({\mathcal G}_F))$ before applying
the transgression map).

Let $E_3=\langle e_1, e_2, e_3\rangle$, with the duals of the $e_i$
corresponding to 
$x_{-1}, x_2, x_5$ respectively. We now define the action on $X_{\mathbb Q_2}$,
using complex coordinates:
\begin{align*}
e_1(z_1, z_2, z_3, z_4, z_5) &= (z_1, z_2, \bar{z}_3, \bar{z}_4,-z_5) \\
e_2(z_1,z_2,z_3, z_4, z_5) &= (-z_1, z_2, -\bar{z}_3, z_4, -\bar{z}_5) \\
e_3(z_1,z_2,z_3,z_4,z_5) &= (z_1,-z_2, z_3, -\bar{z}_4, \bar{z}_5). 
\end{align*}
One notices immediately that the action is free. Moreover if we project onto
the first, second and fifth coordinates, we obtain a fibration
\[
\xymatrix{\Sone\times\Sone \ar[r] & X_{\mathbb Q_2}/E_3 \ar[d] \\
& Y
}
\]
where $Y$ is a 3-dimensional manifold with 
\[
H^*(Y)\cong 
{\mathbb F}_2[e_1,e_2,e_3]/((e_1+e_2)(e_1+e_3), e_2^2,e_3^2).
\]
The fibre generators
transgress to the products $e_1e_2$ and $e_1e_3$.
The cohomology of the orbit space $X_{\mathbb Q_2}/E_3$
can be computed iterating a Gysin sequence argument,
and it has Poincar\'e series equal to
\[
1+3t+6t^2+6t^3+3t^4+t^5.
\]
\end{example}

The example indicates a general fact which we now state.

\begin{theorem}
\label{7.2}
Let $X_F$ be the free $E_n$ space associated to a 
field $F$ which is not formally real,
where  $|\dot F/\dot F^2|=2^n$
and $r=n+{\binom{n}{2}} - \dim H^2(\Gal(F_q/F))$. There is a free
action of $E_n$ on a direct factor $(\Sone)^n\subset X_F$, and
hence we have a fibration 
\[
\xymatrix{(\Sone)^{r-n} \ar[r] & X_F/E_n \ar[d] \\
&(\Sone)^n/E_n
}
\]
\end{theorem}
\begin{proof} 
Let $e_1,\dots ,e_r\in H^1(X_F)$ denote
the 1-dimensional cohomology generators. In the Cartan-Leray
spectral sequence associated to the action of $E_n$, these classes
transgress to 2-dimensional polynomials $\kappa_1,\dots ,\kappa_r
\in H^*(E_n)$. By theorem~\ref{2.3}, 
\[
H^*(F)\cong {\mathcal R}
\cong H^*(E_n)/(\kappa_1,\dots ,\kappa_r) \subset H^*(X_F/E_n)
\]
hence in particular no higher differentials can hit the bottom
edge of the spectral sequence. Notice that $H^*(X_F/E_n)$ is finite,
hence the $\kappa_j$ form a possibly redundant homogeneous system
of parameters in $H^*(E_n)$. We can therefore choose a minimal
subset of parameters in $\{\kappa_1,\dots ,\kappa_r\}$; this subset
must necessarily have $n$ elements in it. By relabelling if necessary,
we can assume that the first $n$ of them form such a subset; note that
they will constitute a regular sequence in $H^*(E_n)$. Now take the
direct factor $Z$ consisting of the first $n$-circles in $X_F$; recalling
that the action we constructed is a product of actions on each of
the circles we see that the factorization of $X_F$ is compatible
with respect to the action of $E_n$. Now consider the action on this
subcomplex; the transgressions form a regular sequence, hence
$H^*_{E_n}(Z)\cong H^*(E_n)/(\kappa_1,\dots ,\kappa_n)$, an algebra
of finite total dimension. Hence the action of $E_n$ is free, and we
obtain the desired conclusions. 
\end{proof}

\end{section}
\begin{section}{Formally Real Fields}
\label{s:formally real fields}

In this section we specialize to fields which are formally real.
As we have seen, the model $X_F$ is endowed with an action of
an elementary abelian group with non-trivial
cyclic isotropy subgroups. We begin our analysis by looking at
a very interesting collection of examples.

Consider the 2-group defined as an extension
\[
1\to ({\mathbb Z}/2)^{\binom{n}{2}}\to T(n)\to ({\mathbb Z}/2)^n\to 1
\]
which can be regarded as the quotient of $W(n)$ obtained by making
the squares of the generators equal to zero. $T(n)$ can also be thought
of as the extension defined by taking $n$ generators of order 2 such that
their commutators are central and also of order 2. More concretely, we
note that $T(2) \cong D_8$, the dihedral group of order 8.
The results of \cite{MS2} can be combined to yield the following
information about fields with \textit{W}-group $T(n)$: 
\begin{proposition}
\label{p:G_F=T(n)}
If ${\mathcal G}_F \cong T(n)$, $F$ is a formally real pythagorean
field with $|\dot{F}/\dot{F}^2| = 2^n$.  
Conversely, the \textit{W}-group of 
such a field is generated by $n$ involutions, and is therefore a quotient
of $T(n)$.
\end{proposition}
\begin{proof}
By \cite[2.11]{MS2}, $F$ is pythagorean if and only if ${\mathcal
G}_F$ is generated by involutions, and $T(n)$ is, by definition, generated by
involutions.  By \cite[2.7.2]{MS2}, $F$ is formally real if 
there are involutions of ${\mathcal G}_F$ which do not lie in the
Frattini subgroup.  Again, for $T(n)$ this condition obviously holds. 
The condition on $|\dot{F}/\dot{F}^2|$ follows immediately from the
fact that ${\mathcal G}_F/\Phi({\mathcal G}_F)$ is dual to
$|\dot{F}/\dot{F}^2|$. The partial converse follows from the fact that
\emph{any} \textit{W}-group generated by $n$ involutions is a quotient of
$T(n)$ and \cite[2.11]{MS2}.
\end{proof}

This statement can be strengthened to a characterization
of fields with $T(n)$ as \textit{W}-group. The proof of this
characterization involves some field theory which may not be totally
familiar, but most of the results which we use can 
be found in \cite{Lam2}.  For the reader's convenience we recall
briefly the relevant definitions.  A subset $P \neq F$ of a field $F$
is said to be an \textsl{ordering} if $P$ is closed under addition and
multiplication and $P \cup -P = F$. Elements of $P$ are said to be
\textsl{positive} with respect to the ordering. The motivating example
is $F = {\mathbb R}$ and $P = \{x \geq 0\}$; by
the Artin-Schreier theorem \cite[1.5]{Lam2}, an ordering of $F$ exists
if and only if $F$ is formally real. A \textsl{preordering}
\cite[p.~2]{Lam2} weakens the requirement $P \cup -P = F$ to simply
$P \supset F^2$; it follows that every preordering contains $\Sigma
F^2$, the set of all sums of squares of elements of $F$, and this set
is called the \textsl{weak preordering} of $F$.  Roughly speaking, a
preordering $T$ of $F$ has the \textsl{strong approximation property}
(abbreviated SAP) if, given disjoint subsets of the set of orderings
containing $T$, it is possible to choose an element of $F$ positive
with respect to the first set of orderings and negative with respect
to the second. $F$ itself is said to be \textsl{SAP} if the weak
preordering of $F$ has the SAP \cite[p.~126]{Lam2}. 
It turns out that a formally real pythagorean field with
$|\dot F/\dot F^2|=2^n$ is SAP if and only if $F$ admits exactly $n$
orderings \cite[17.4]{Lam2}.  This last is a key fact we will use in
proving the characterization below:

\begin{theorem}
\label{t:G_F=T(n)}
${\mathcal G}_F \cong T(n)$ if and only if $F$ is a formally real
pythagorean SAP field with $|\dot F/\dot F^2|=2^n$.  Moreover,
for each $n\in\mathbb N$ there exists such a field.
\end{theorem}
\begin{proof} 
The existence of formally real pythagorean SAP fields with
$|\dot F/\dot F^2|=2^n$ for each $n\in\mathbb N$ was proved in
\cite{Bro}, Satz (4).  We therefore turn to the proof of the ``if and
only if'' part.
\begin{proof}[Necessity]
Assume that ${\mathcal G}_F \cong T(n)$.  By
proposition~\ref{p:G_F=T(n)}, we need only show that $F$ is SAP, i.e.\
that $F$ has exactly $n$ orderings.  By \cite[2.10]{MS2}, there is a 
bijection between the set of all orderings of $F$ and the non-trivial
cosets $\sigma \Phi({\mathcal G}_F)$, where $\sigma$ is an involution.
(By ``non-trivial'' we mean $\sigma \notin \Phi({\mathcal G}_F)$.)
Since ${\mathcal G}_F \cong T(n)$, let us write the involutions
generating ${\mathcal G}_F$ as $\sigma_1,\ldots ,\sigma_n$.  We claim that the
only non-trivial cosets $\sigma \Phi({\mathcal G}_F)$ with $\sigma$ an
involution are those with $\sigma = \sigma_i$ for some $i$, or
equivalently, that the maximal elementary abelian subgroups of
${\mathcal G}_F$ are the subgroups $\langle \sigma_i, \Phi({\mathcal
G}_F) \rangle$.  To see 
this, write the general element $\sigma \in {\mathcal G}_F$ as a
product:
\[
\sigma = \sigma_1^{\epsilon_1}\sigma_2^{\epsilon_2}\cdots
\sigma_n^{\epsilon_n}\prod [\sigma_i,\sigma_j]^{\epsilon_{ij}}
\]
where $\epsilon_i,\epsilon_{ij}\in\{0,1\}$, and note that if $\sigma
\notin {\mathcal G}_F$, $\epsilon_i = 1$ for at least one $i$,
while if $\epsilon_i = 1$ for more than one $i$, $\sigma^2 \neq 1$. 
\renewcommand{\qedsymbol}{}
\end{proof}
\begin{proof}[Sufficiency]
The converse is more interesting. Let $F$ be a formally real
pythagorean SAP field, with $|\dot F/\dot F^2|=2^n$ and suppose that
${\mathcal G}_F\ne T(n)$. 
Then we want to show that there exist more than $n$ orderings of $F$. 
This will require some field theory; to start choose
$a_1,\dots ,a_n\in \dot F$ such that 
$\sigma_i(\sqrt{a_j})/\sqrt{a_j}=-1$ if $i=j$ and is equal to $1$
if $i\ne j$. Then $a_i\in \cap_{i\ne j} P_j-P_i$, where
$P_i=\{f \mid \sqrt{f}^{\sigma_i}=\sqrt{f}\}$ are orderings corresponding
to the $\sigma_i$ (see \cite[p.~522]{MS2}). (The 
$a_i$ exist since  $F$ is assumed to be SAP.)
It follows that the cosets $[a_1],\dots ,
[a_n]$ form a basis for the vector space $\dot F/\dot F^2$.
If ${\mathcal G}_F\ne T(n)$, then ${\mathcal G}_F$ is a proper quotient
of $T(n)$ (\ref{p:G_F=T(n)}), and there exists a nontrivial relation
between commutators 
$[\sigma_i,\sigma_j]$, $1\le i<j\le n$. Without loss of generality
we may assume that our relation $\gamma = [\sigma_1,\sigma_2]\gamma_1\in T(n)$
is a nontrivial product of commutators of the form 
$[\sigma_i,\sigma_j]$, $i<j$, and that $[\sigma_1,\sigma_2]$ does not
enter the expression for $\gamma_1$. Since $\gamma$ is a relation, the image
$\bar\gamma$ of $\gamma$ in ${\mathcal G}_F$ is $1$. From 
\cite[2.21]{MS1} we can conclude that the quaternion algebra
\[
A=\langle \alpha_1,\alpha_2 \mid \alpha_1^2=a_1, \alpha_2^2=a_2,
\alpha_1\alpha_2=-\alpha_2\alpha_1\rangle
\]
does not split. Applying the standard criterion for the splitting
of quaternion algebras over $F$ (see \cite{Lam1}), we see that
$a_2\notin D\langle 1,-a_1\rangle$, where by $D\langle 1,-a_1\rangle$
we mean the set of 
values of the quadratic form $x^2-a_1y^2$ over $F$.
However, using \cite[1.6]{Lam2}, we see that 
$D\langle 1,-a_1\rangle$ is the intersection of all orderings $P$ of
$F$ such that $-a_1\in P$. If the set of all orderings were just
$P_1,\ldots ,P_n$ (the ones determined by the cosets $\sigma_i\Phi
({\mathcal G}_F)$), then we would have that this intersection is equal
to $P_1$, but $a_2\in P_1$, a contradiction. Therefore we see that
there exists some ordering $P$ of $F$ such that $P\ne P_1,\ldots
,P_n$.   
\renewcommand{\qedsymbol}{}
\end{proof}
This finishes the
proof of the fact that ${\mathcal G}_F\cong T(n)$ 
if and only if $F$ is a formally real pythagorean SAP field with
$|\dot F/\dot F^2|=2^n$.
\end{proof}

\begin{remark} 
The theorem above also follows from \cite[prop.~4.1]{CS}, where
the results in \cite{MSm1} were applied, but we have chosen to provide a
direct proof.
\end{remark}

We now study the cohomology of the groups $T(n)$
using the topological models previously described in \ref{t:X_F}.
If $H^*(({\mathbb Z}/2)^n)= {\mathbb F}_2[x_1,\dots ,x_n]$, then the $k$-invariants
for the defining extension for $T(n)$
are precisely $x_ix_j$, with $i<j$. We consider an
action of $E_n=({\mathbb Z}/2)^n$ on $Z_n=(\Sone)^{\binom{n}{2}}$ defined
as follows:
if $\langle e_1,\dots ,e_n\rangle$ is a basis for $E_n$, then
we define 
\[
e_k(z_{ij})=
	\begin{cases}
		\bar{z}_{ij},& \text{if $k=i$;} \\
                -\bar{z}_{ij},& \text{if $k=j$;} \\
                z_{ij},& \text{otherwise.}
	\end{cases}
\]

In the notation of \ref{t:X_F}, if $F$ is a field with ${\mathcal
G}_F\cong T(n)$, then we can assume that $X_F=Z_n$, with the above
action of $({\mathbb Z}/2)^n$.  We now analyze the equivariant
cohomology ring $H^*_{E_n}(X_F)$.  To do this we need to
understand the ``singular set'' of the action, i.e.\ 
the subspace of elements in $Z_n$ which are not freely permuted by
the group $E_n$. In sufficiently high dimensions the equivariant cohomology
is completely determined on this subspace (\ref{ec:large-dimensions}). 

Let $C_k=\langle e_k\rangle\subset E_n$; we
examine its fixed-point set. 
Fixing $k$, there are $k-1$ coordinates of the form $z_{ik}$ and
$n-k$ coordinates of the form $z_{kj}$ in $Z_n$. The group acts on each
coordinate independently. Hence an ${\binom{n}{2}}$-tuple $z_{ij}$ will be
fixed by $C_k$ if and only if $z_{ik}\in \{ i,-i\}$ for $i<k$ and
$z_{kj}\in \{1,-1\}$ for $j>k$ (note that the action on the other 
coordinates is trivial). Hence we obtain
\[
Z_n^{C_k}\cong \{1,-1\}^{n-k}\times \{i,-i\}^{k-1}\times 
(\Sone)^{\binom{n-1}{2}}.
\]
Now let $E_{n,k}=E_n/C_k$; then this group acts on $Z_n^{C_k}$, and we
can express this space as a $E_{n,k}$ space as follows:
\[
Z_n^{C_k}\cong E_{n,k}\times Z_{n-1}
\]
where $E_{n,k}\cong E_{n-1}$ acts as before on $Z_{n-1}$ (here we have all
the coordinates $z_{ij}$ which do not involve $k$; the basis we use for
$E_{n,k}$ is simply the one obtained by omitting $e_k$). We need some 
understanding of the geometry of these fixed point sets. We have

\begin{proposition}
\label{7.4}
$Z_n^{C_k}\cap Z_n^{C_l}=\varnothing$ if $k\ne l$ and any cyclic
subgroup
in $E_n$ distinct from $C_1,\dots, C_n$ acts freely on $Z_n$.
\end{proposition}

\begin{proof} Assume that $k<l$ and consider the coordinate
$z_{kl}$; then the fixed point set for $C_k$ has $1$ or $-1$ in this
coordinate and the fixed point set for $C_l$ has $i$ or $-i$ in this
coordinate. Hence they cannot intersect. Now let $U\subset E_n$ be any
cyclic subgroup distinct from the $C_i$. It will contain an element
of the form $x_{i_1}x_{i_2}\dots x_{i_r}$ where $r>1$ and 
$i_1<i_2<\dots <i_r$. This will act on the coordinate $z_{i_1i_2}$
as multiplication by $-1$, hence it acts freely on all of $Z_n$.
\end{proof} 

We have shown that the singular set ${\mathcal S}$ of the action is 
a disjoint union of precisely $n$ subspaces, which, as $E_n$
spaces, are of the
form $E_{n,k}\times Z_{n-1}$. Since 
$H^*_{E_n}(E_{n,k}\times Z_{n-1})\cong H^*_{C_k}(Z_{n-1})$
(\ref{ec:frobenius}), we obtain 
\[
H^*_{E_n}({\mathcal S})\cong H^*_{C_1}(Z_{n-1})\oplus\dots 
\oplus H^*_{C_n}(Z_{n-1}).
\]
Each of these factors has Poincar\'e series equal to
$(1+t)^{\binom{n-1}{2}}/(1-t)$; hence we deduce that the Poincar\'e
series for the equivariant cohomology of the singular set is
$n(1+t)^{\binom {n-1}{2}}/(1-t)$. In dimensions larger than
${\binom{n}{2}}$ this will coincide with the equivariant cohomology
(\ref{ec:large-dimensions}). Using this we obtain that there is a
polynomial $s_n(t)$ of degree ${\binom{n}{2}}$ such that the
cohomology of $T(n)$ is given by the expression
\begin{equation}
\label{e:PS(T(n))}
\frac{(1-t)s_n(t) + n2^{\binom{n-1}{2}}t^{{\binom{n}{2}}+1}}
{(1-t^2)^{\binom{n}{2}}(1-t)}.
\end{equation}

The polynomial $s_n(t)$ is simply the Poincar\'e series for the
equivariant cohomology $H^*_{E_n}(Z_n)$ through dimension ${\binom{n}{2}}$.
For example, if $n=2$ this is $1+2t$ and we obtain the series
\[
\frac{(1-t)(1+2t) + 2t^2}{(1-t^2)(1-t)}=\frac{1}{(1-t)^2}
\]
which agrees with that for the cohomology of the dihedral group
$D_8\cong T(2)$.

It is worthwhile to observe that there is a one dimensional class,
$\beta = x_1+x_2+\dots +x_n\in H^1(T(n))\cong H^1(({\mathbb Z}/2)^n)$
which restricts non trivially on $H^1_{C_i}(Z_{n-1})$ for all
$i=1,\dots ,n$. On each of these summands this class restricts to
a polynomial class coming from the cohomology of the cyclic isotropy
subgroup. In sufficiently high dimensions, multiplication by this
class induces a periodicity isomorphism. Hence if we take 
$\{\mu_{ij}\}$ and adjoin $\beta$, the quotient
$H^*(T(n))/(\mu_{ij},\beta)$ is finite dimensional. Hence we conclude
that $\beta$ together with the $\mu_{ij}$, $i<j$ form a 
\textsl{homogeneous
system of parameters} for $H^*(T(n))$. We also record the following
simple consequence of our calculations:

\begin{corollary}
\label{7.6}
Let $F$ be a field such that ${\mathcal G}_F=T(n)$; then
\[
H^*(F)\cong {\mathbb F}_2[x_1,\dots ,x_n]/(x_ix_j)\cong
{\mathbb F}_2[x_1]\bar{\oplus}\cdots\bar{\oplus}
{\mathbb F}_2[x_n],
\]
where $\bar{\oplus}$ indicates that we identify the unit elements in
the rings involved in the direct sum. In particular, multiplication by $\beta$
is an isomorphism. 
\end{corollary}

To obtain further calculations, one needs to understand the series $s_n(t)$.
In some instances this can be obtained by direct geometric methods. Consider
the case $n=3$. $Z_3$ is a 3-torus with an action of $E_3$ such that the
singular set is homeomorphic to the disjoint union of 12 circles. One can readily
see that there is a CW-decomposition for $Z_3$ which is permuted by this 
action (this is the product of the natural cellular structure for each
circle). There are $64$ three-dimensional cells and $192$ 
two-dimensional cells,
all of which are freely permuted by $E_3$. In fact the cellular chain complex
has the following form:

\[{\mathbb F}_2[E_3]^8\to 
{\mathbb F}_2[E_3]^{24}\to [\oplus_{i=1}^3{\mathbb F}_2[E_3/C_i]^2]
\oplus {\mathbb F}_2[E_3]^{18}\to~~~~
\]
\[
\hskip 1in \to [\oplus_{i=1}^3
{\mathbb F}_2[E_3/C_i]^2]\oplus{\mathbb F}_2[E_3]^2.
\]

The singular set is of codimension 2, and the pair $(Z_3,{\mathcal S})$ 
is a relative 3-manifold; 
hence we deduce that $H^3(Z_3/E_3)\cong H^3_G(Z_3,{\mathcal S})\cong {\mathbb F}_2$. 
Using this fact and the long exact sequence associated to the equivariant pair
$(Z_3,{\mathcal S})$, 
one can show that 
$s_3(t)=1+3t+5t^2+6t^3$. Note that the coefficients 
arise from the 
isomorphisms
$H^1_{E_3}(Z_3)\cong H^1(E_3)$, $H^2_{E_3}({\mathcal S})\cong H^2_{E_3}(Z_3)
\oplus H^3(Z_3/E_3)$
and $H^3_{E_3}(Z_3)\cong H^3_{E_3}({\mathcal S})$.
Substituting in our formula above, we obtain
that the Poincar\'e series for the cohomology of $T(3)$ is 
\[
\frac{{(1-t)(1+3t+5t^2+6t^3) +6t^4}}
{(1-t)(1-t^2)^3}=\frac{1+2t+2t^2+t^3}{(1-t)(1-t^2)^3}.
\]
The cohomology of this group has already been computed (see \cite{C1},
group number 144) and it agrees with the above;
it is detected on the 3 conjugacy
classes of maximal elementary abelian subgroups.

For later use we point out that in proving \ref{t:G_F=T(n)} we established
that $T(n)$ has precisely
$n$ maximal 2-tori, all of rank ${\binom{n}{2}}+1$. 
If $\sigma_1,\dots ,\sigma_n$ generate the group, then the tori are of the
form
\[
U_k=\langle [\sigma_i,\sigma_j], \sigma_k \mid 1\leq i < j \leq n \rangle.
\]
Note that they all intersect along the central
elementary abelian subgroup generated by the commutators, and they
are non-conjugate.
These groups are all self-centralizing, and their Weyl groups are
isomorphic to $({\mathbb Z}/2)^{n-1}$.

To have a better understanding of the situation 
we will describe a rather interesting relationship between
the groups
$T(n)$ and $W(n-1)$. If $T(n)$ is generated by $\sigma_1,\dots ,\sigma_n$,
consider the subgroup $P(n)$ generated by the elements $\sigma_1\sigma_2,\dots ,
\sigma_1\sigma_n$ and the ${\binom{n}{2}}$ central commutators of order two.
This is an index two subgroup, which we claim is isomorphic to $W(n-1)$.

\begin{proposition}
\label{7.7}
The group $T(n)$ is isomorphic to the semi-direct product
$W(n-1)\times_T{\mathbb Z}/2$, where the element of order two acts by
inverting the $n-1$ generators in $W(n-1)$.
\end{proposition}
\begin{proof} 
We will use the group $P(n)$ to establish this
isomorphism. Note that $T(n)/P(n)$ is generated by the class represented
by $\bar{\sigma}_1$; this of course splits as $\sigma_1$ is an
element of order two in $T(n)$. Hence $T(n)$ is a semi-direct product, 
and it will suffice to show that $P(n)\cong W(n-1)$ and identify
the action of $\sigma_1$. 

Observe that $W(n-1)$ is the universal group in the category of all
2-groups with $n-1$ generators and satisfying the following two conditions:
(1)~all squares are central, (2)~all elements have order $\le 4$. Given that 
$P(n)\subset T(n)$, it will automatically satisfy both of these conditions
and hence is a homomorphic image of $W(n-1)$. However by comparing orders
we see that this in fact they must be isomorphic. 
The action
is determined by conjugation with $\sigma_1$. Note that
$\sigma_1\sigma_1\sigma_i\sigma_1=\sigma_i\sigma_1=(\sigma_1\sigma_i)^{-1}$,
whence the proof is complete.
\end{proof} 

We would like to understand \ref{7.7} from the point of view of field theory.
Assume that $F$ is a field such that ${\mathcal G}_F=T(n)$. We claim that $P(n)$
can be interpreted as $\Gal(F^{(3)}/F(\sqrt{-1}))$. We have already
remarked that under these conditions $F$ is a formally real, pythagorean
SAP field. In particular $\sqrt{-1}\notin \dot F$. From \cite[p.~521]{MS2},
we see that for each $\sigma_i$, $\sigma_i(\sqrt{-1})=-\sqrt{-1}$. Hence
for each $\sigma_i,\sigma_j$ we have $\sigma_i\sigma_j (\sqrt{-1})
=\sqrt{-1}$, and therefore 
$P(n)\subset \Gal(F^{(3)}/F(\sqrt{-1}))$; as it has index two
in ${\mathcal G}_F$ this must be an equality. More explicitly, we claim that
$P(n)$ can be identified with the \textit{W}-group associated to $F(\sqrt{-1})$.
Indeed, from \cite[thm.~3.4, p.~202]{Lam1} we have an exact sequence
induced by the inclusion of $F$ in $F(\sqrt{-1})$:
\[
1\to \{\dot F^2,-\dot F^2\}\to \dot F/\dot F^2
\buildrel{\epsilon}\over\longrightarrow 
{\dot F}(\sqrt{-1})/
{\dot F}(\sqrt{-1})^2\xrightarrow{N}\dot F/\dot F^2
\]
where $\epsilon$ is the map induced by inclusion, and $N$ is the
homomorphism induced by the norm. In our case $N$ is trivial, and
$\epsilon$ is surjective. Note that $|\dot F(\sqrt{-1})/\dot
F(\sqrt{-1})^2|=2^{n-1}$, hence ${\mathcal G}_{F(\sqrt{-1})}$ is a 
homomorphic image of $W(n-1)$. On the other hand we see that
$F(\sqrt{-1})^{(2)}$ (the compositum of all quadratic extensions
of $F(\sqrt{-1})$) is just $F^{(2)}(\sqrt{-1})$. Since $F^{(3)}$ is
Galois over $F(\sqrt{-1})$, we see that $F^{(3)}\subset 
F(\sqrt{-1})^{(3)}$. Hence $P(n)=\Gal(F^{(3)}/
F(\sqrt{-1}))$ is a homomorphic image of ${\mathcal G}_{F(\sqrt{-1})}$.
We conclude that $P(n)\cong {\mathcal G}_{F(\sqrt{-1})}\cong W(n-1)$.
We have identified the class of the extension
associated to $P(n)\subset T(n)$: it is the class $\alpha
\in H^1(T(n))$ which corresponds to the element
$[-1]\in \dot F/\dot F^2\cong H^1(T(n))$. Geometrically we see
that $\alpha =\beta$, a class which ``lives'' on all components of
the singular set of the action on $X_F$.

Associated to the extension class $\alpha$ we have a Gysin sequence
\[
\dots H^i(T(n))\xrightarrow{\cup\alpha}
H^{i+1}(T(n))\xrightarrow{\Res} H^{i+i}(P(n)) \xrightarrow{\Tr}
H^{i+1}(T(n))\to\dots
\]
If multiplication by $\alpha$ were injective, then this sequence
would collapse, i.e.\ the cohomology of $P(n)$ would be totally
nonhomologous to zero, and the Poincar\'e series for the
cohomology of $T(n)$ would be of the form
$p_{n-1}(t)(1-t)^{-1}$
where $p_{n-1}(t)$ is the Poincar\'e series for $W(n-1)$
(\ref{ec:TNHZ}). Moreover 
the homogeneous system of parameters formed by
adjoining the class $\alpha$,
to the regular sequence $\{\mu_F\}$ would be regular, and 
hence $H^*(T(n))$
would be Cohen-Macaulay. Conversely if $H^*(T(n))$ is Cohen-Macaulay,
any homogeneous system of parameters is a regular sequence, and
this in turn implies that multiplication by $\alpha$ is injective.
Moreover if $H^*(T(n))$ is Cohen-Macaulay, it must be detected
on the centralizers of elementary abelian subgroups of maximal
rank (a result appearing in \cite{C2}); in this case we obtain detection
on elementary abelian subgroups.
Conversely, note that the homogeneous
system of parameters consisting of $\alpha$ and the $\mu_F$ 
restricts to regular sequences on all the maximal elementary abelian
subgroups; hence if they detect the entire cohomology this must
be a regular sequence.

 From our preceding calculations we see that in the Gysin sequence
for the extension 

$$1\to P(3)\to T(3)\to {\mathbb Z}/2\to 1$$
the cohomology of $P(3)$ is totally non-homologous to zero, and hence
it collapses to yield the Poincar\'e series described previously. We also
have explicit detection on elementary abelian subgroups.
We claim that this simple scheme breaks down for $T(4)$
and from there on. Assume that multiplication by $\alpha$ is
injective; then the Poincar\'e series for $H^*(T(4))$ would
be given by 
\[
\frac{1+3t+8t^2+12t^3+8t^4+3t^5+t^6}{(1-t)(1-t^2)^6}.
\]
In particular if we divide out by the regular sequence
$\mu_F$, we see that in dimension 7 there must be precisely
$1+3+8+12+8+3+1=36$ linearly independent elements. However, from our
previous calculation of Poincar\'e series,
there can only be 32.  This
contradiction means that multiplication 
by $\alpha$ cannot be injective, and so we deduce that $H^*(T(4))$
is not Cohen-Macaulay. 

More generally, there is a natural
inclusion of pairs $(T(n),P(n))\subset (T(n+1),P(n+1))$ which gives
rise to a commutative diagram of restriction maps:
\[
\begin{CD}
H^*(T(n+1)) 	@>>> 	H^*(T(n)) 	\\
@VVV 			@VVV 		\\
H^*(P(n+1)) 	@>>> 	H^*(P(n))	\\
\end{CD}
\]
The two horizontal arrows are naturally split
and hence are automatically surjective.
Now injectivity of multiplication by $\alpha$ is 
equivalent to surjectivity
of the restriction map $H^*(T(n))\to H^*(P(n)).$
Hence we deduce that if multiplication by $\alpha$
is injective for $T(n+1)$, then it must be true for $T(n)$
(e.g.\ $T(3)$ and $T(2)$). Hence as it fails for $T(4)$, 
it must fail for $T(n)$ for all $n\ge 4$. We therefore
have
\begin{proposition}
\label{7.8}
For $n\ge 4$, $H^*(T(n))$ has depth equal to one less
than the rank of the group. Hence it is not Cohen-Macaulay and
it is not detected on elementary abelian subgroups.
\end{proposition}

Using the insight acquired from our analysis of $T(n)$,
we now make explicit the role played by the class $[-1]$
in the cohomology of ${\mathcal G}_F$ for any formally real field.

\begin{proposition}
\label{7.9}
Let $F$ be formally real. Then the class $[-1]\in H^1({\mathcal G}_F)$
restricts non-trivially on any cyclic subgroup of the form
$\langle \sigma\rangle$, where $\sigma\in {\mathcal G}_F-\Phi({\mathcal G}_F)$ is a 
non-central involution. Furthermore if $\alpha\in H^*({\mathcal G}_F)$ 
satisfies $[-1]\cup\alpha=0$, then $\alpha$ is nilpotent.
\end{proposition}
\begin{proof}
Let $\sigma\in {\mathcal G}_F$ denote 
a non-central involution. Now let $P$ be the ordering corresponding 
to $\sigma$, i.e.\ $P=\{p\in \dot F \mid \sqrt{p}^{\sigma}=\sqrt{p}\}$. Since
$-1\notin P$, we see that $\sqrt{-1}^{\sigma}=-\sqrt{-1}$, and therefore
$[-1](\sigma)\ne 0$. This shows that $[-1]$ restricts non-trivially
on any such cyclic subgroup $\langle\sigma\rangle$ (note that in particular
this shows that $[-1]$ is not nilpotent). The maximal elementary
abelian subgroups in ${\mathcal G}_F$ are of the form 
$E_{\sigma}=\langle\Phi({\mathcal G}_F),\sigma\rangle$, where is $\sigma$ is a non-central
involution as before. Hence in particular we have shown that $[-1]$
restricts non-trivially to every maximal elementary abelian subgroup.
Thus, if $[-1]\cup\alpha=0$, we know that $\alpha$ must restrict 
trivially on all these subgroups, and by Quillen's theorem
(see \cite{AM}) must be nilpotent.
\end{proof}

We are interested in investigating under what conditions $[-1]$ is
not a zero divisor. We have

\begin{lemma}
\label{7.10}
Let $F$ be a formally real field. Then, if $[-1]$ is not
a zero divisor in $H^*({\mathcal G}_F)$, $F$ must be pythagorean.
\end{lemma}
\begin{proof}
Suppose that $F$ is a formally real, non-pythagorean field.
Then there exists an element $a\in\dot F$ such that $a\notin\dot F^2$
but is a sum of two (non-zero) squares in $F$, say $a=x^2 +y^2$. If we rewrite
this as $a-x^2-y^2=0$, then we obtain the relation
$[a][-x^2][-y^2]=0$ 
in $H^*({\mathcal G}_F)$. Here we are using a well-known identity
in Milnor \textit{K}-theory (see \cite[p.~320]{M}) and the fact that we have
identified it with the subring ${\mathcal R}\subset H^*({\mathcal G}_F)$
generated by 1-dimensional classes (see \ref{2.3}). Expanding this
out by using the \textit{K}-theoretic relations, we obtain

$$0=[a][-x^2][-y^2]=[a]([-1]+[x^2])([-1]+[y^2])=[a][-1]^2$$
from which we deduce that $[-1]$ is necessarily a zero divisor,
a contradiction. Hence $F$ must be pythagorean.
\end{proof}

Now given any formally real field $F$ we can construct an index 2
subgroup ${\mathcal H}\subset {\mathcal G}_F$ corresponding to
$[-1]\in\dot F/\dot F^2$. We can use this subgroup 
to obtain a simple condition that implies the
Cohen-Macaulay property.

\begin{proposition}
\label{7.11}
If $F$ is formally real and $[-1]$ is not a zero
divisor, then $H^*({\mathcal G}_F)$ is
Cohen-Macaulay.
\end{proposition}
\begin{proof} 
Assume that $[-1]$ is not a zero divisor.
Then using the Gysin sequence we can express
$$H^*({\mathcal H})\cong H^*({\mathcal G}_F)/([-1]).$$
Now recall that the sequence
$\{\mu_F\}$ in $H^*({\mathcal G}_F)$ is regular; but more is true: 
in fact the classes
$\mu_F$ restrict non-trivially to $\Phi({\mathcal G}_F)$, the maximal
central elementary abelian subgroup. The index 2 subgroup 
${\mathcal H}$ contains this central elementary abelian
subgroup, and hence the $\mu_F$ restrict to a regular sequence
in its cohomology. This means that $\{\mu_F\}$ is a regular sequence
in the quotient $H^*({\mathcal G}_F)/([-1])$, and we conclude that
$\{\mu_F,[-1]\}$ is a regular sequence in $H^*({\mathcal G}_F)$, and so
it is Cohen-Macaulay.
\end{proof} 

Next we observe that the centralizers
of the maximal elementary abelian 
subgroups  $E_{\sigma}\subset {\mathcal G}_F$ are simply the centralizers of the 
non-central involutions $\sigma\in {\mathcal G}_F-\Phi({\mathcal G}_F)$.
An immediate consequence of lemma 2.1 in \cite{MMS} is that these
centralizers are precisely the elementary abelian subgroups themselves
i.e.\ they are self-centralizing.
Applying the detection results in \cite{C2} we obtain a simple
consequence of \ref{7.11}. 

\begin{proposition}
\label{7.12}
If $F$ is formally real then $[-1]$ is not a zero divisor
if and only if $H^*({\mathcal G}_F)$ 
is detected on the 
cohomology of its elementary abelian subgroups.
\end{proposition}
\begin{proof}
According to \cite{C2}, if $H^*({\mathcal G}_F)$
is Cohen-Macaulay, it is detected on centralizers of maximal 2-tori.
On the other hand, if we assume detection on elementary abelian subgroups,
then $[-1]$ cannot be a zero divisor. Indeed, we have shown that it
restricts non-trivially on all the maximal elementary abelian subgroups,
hence given any $x\in H^*({\mathcal G}_F)$, there must exist an $E_{\sigma}$
on which $[-1]\cup x$ restricts non-trivially, as this is true for $x$.
\end{proof}

It remains to identify the subgroup ${\mathcal H}$ in terms of
field-theoretic data. Following the arguments given after \ref{7.7}, one
can show that ${\mathcal H}$ is a quotient of the \textit{W}-group 
${\mathcal G}_{F(\sqrt{-1})}$. 
Indeed the arguments provided there can be extended to provide
the following entirely analogous characterization of fields where
these two groups agree (see \cite{MSm2} for details).

\begin{theorem}
\label{7.13}
If $F$ is not formally real then the index two subgroup
${\mathcal H}\subset {\mathcal G}_F$ corresponding to $[-1]$
is isomorphic to ${\mathcal G}_{F(\sqrt{-1})}$ if and only if
$F$ is pythagorean. Furthermore, in that case
the group ${\mathcal G}_F$
can be expressed as a semi-direct product
$${\mathcal G}_F\cong {\mathcal G}_{F(\sqrt{-1})}\times_T{\mathbb Z}/2$$
where the involution acts by inverting a suitable collection
of generators for ${\mathcal G}_{F(\sqrt{-1})}$.
\end{theorem}

Now let $F$ denote a pythagorean field.
Evidently $F(\sqrt{-1})$ is not formally
real, hence its \textit{W}-group satisfies the 2C property. This has an 
interesting geometric interpretation. Take the model $X_F$ with an
action of $E_n$. Then there exists a hyperplane $H_{-1}\subset E_n$
such that $X_F$ with the action of $H_{-1}$ is a model for 
${\mathcal G}_{F(\sqrt{-1})}$ (indeed, we can identify this
kernel with 
$H^1({\mathcal G}_{F(\sqrt{-1})})\subset H^1({\mathcal G}_F)\cong E_n$).
Note that $\Phi ({\mathcal G}_F)=\Phi ({\mathcal G}_{F(\sqrt{-1})})$.
By our previous observation, $H_{-1}$ must act freely on $X_F$.
Let ${\mathcal I}=\{C_1,\dots ,C_t\}$ denote the set of all isotropy
subgroups for the action of $E_n$ on $X_F$. Then $H_{-1}$ is a 
hyperplane which does not contain any of these subgroups.

Putting things together we obtain:

\begin{theorem}
\label{7.14}
If $F$ is a formally real field then the following conditions are
equivalent 
\begin{enumerate}

\item $F$ is pythagorean and $H^*({\mathcal G}_F)$ is Cohen-Macaulay.

\item $F$ is pythagorean and $H^*({\mathcal G}_F)$ is detected on
its elementary abelian subgroups.

\item $[-1]\in H^1({\mathcal G}_F)$ is not a zero divisor.
\end{enumerate}
\end{theorem}
\begin{proof} 
By our previous results it suffices
to show that (1) implies (3).
If we have the Cohen-Macaulay
condition, any non-redundant homogeneous system of parameters forms a regular
sequence. We claim that this is true for $\{\mu_F,[-1]\}$. It suffices
to show that multiplication by $[-1]$ on $H^*({\mathcal G}_F)/(\mu_F)$
is an isomorphism in sufficiently large dimensions. Consider the
action of $E_n$ on $X_F$; we can identify the quotient above with the
equivariant cohomology $H^*_{E_n}(X_F)$. Similarly we can identify
$H^*({\mathcal G}_{F(\sqrt{-1})}/(\mu_F)$ with $H^*(X_F/H_{-1})$. The
Gysin sequence is recovered on equivariant cohomology as the following
long exact sequence:
\[
\cdots
H^i(X_F/H_{-1})\rightarrow
H^i_{E_n}(X_F)\xrightarrow{\cup [-1]}
H^{i+1}_{E_n}(X_F)\rightarrow H^{i+1}(X_F/H_{-1})\cdots
\]
We observe that $H^i(X_F/H_{-1})=0$ for $i>\dim X_F$, whence
multiplication by $[-1]$ is eventually injective on the quotient,
and
hence we can conclude that $[-1]$ is not a zero-divisor.
\end{proof}

We will now consider examples of fields where the above conditions
are satisfied. We first recall from \cite[p.~45]{Lam2} that a field
$F$ is said to be \textsl{superpythagorean} if it satisfies the 
following conditions. $F$ is a formally real field
with the property that for any
set $S$ containing $\dot F^2$ but such that $-1\notin S$, if
$S$ is a subgroup of index 2 in $\dot F$, then $S$ is an ordering
on $F$. It is easy to see that superpythagorean fields are pythagorean
(we refer the reader to \cite{Lam2}, appendix to section 5, for additional
information and references).
Nice examples of such fields are given by
$F={\mathbb R} ((t_1))((t_2)) \dots ((t_n))$, the field of iterated power
series over ${\mathbb R}$. In this case $|\dot F/\dot F^2|=2^{n+1}$, with 
a basis given by the classes $[-1], [t_1],\dots , [t_n]$.

We will compute the cohomology of any \textit{W}-group arising from such
a field. For this we make use of the following group-theoretic
characterization described in \cite{MS2}. If $F$ is superpythagorean with
$|\dot F/\dot F^2|=2^n$, then 
\[
{\mathcal G}_F\cong ({\mathbb Z}/4)^{n-1}\times_T{\mathbb Z}/2
\]
a semidirect product where the element of order two acts by inverting
the $n-1$ standard generators in $({\mathbb Z}/4)^{n-1}$. We shall denote
this group by $S(n)$; then the Frattini subgroup is generated by the
elements of order 2 in $({\mathbb Z}/4)^{n-1}$ and we can express it as
a central extension

$$1\to ({\mathbb Z}/2)^{n-1}\to S(n)\to ({\mathbb Z}/2)^n\to 1.$$
If $x_1,\dots ,x_n$ generate the cohomology of the quotient group,
then we can assume that the $k$-invariants of this extension are given by
$$x_2(x_2+x_1), x_3(x_3+x_1),\dots , x_n(x_n+x_1).$$
Note that these elements form a regular sequence in ${\mathbb F}_2[x_1,\dots ,x_n]$.
Now let $X_F$ denote the space
constructed as before; in this case $X_F\simeq (\Sone)^{n-1}$
with an action of $E_n=({\mathbb Z}/2)^n$. Consider the spectral sequence
for computing the equivariant cohomology of $X_F$; it has $E_2$ term

$${\mathbb F}_2[x_1,\dots ,x_n]\otimes \Lambda (u_1,\dots ,u_{n-1})$$
where the exterior generators on the fiber transgress to
the $k$-invariants above. As these form a regular sequence, the spectral
sequence collapses at $E_3$, and furthermore the edge homomorphism 
induces an isomorphism of algebras
$$E_3^{*,0}\cong H^*_{E_n}(X_F)$$
(see \cite{Ca} or \cite{D}). In this case we simply obtain

$$H^*_{E_n}(X_F)\cong H^*(F)\cong 
{\mathbb F}_2[x_1,\dots ,x_n]/x_2(x_1+x_2),\dots ,x_n(x_1+x_n).
$$
Note that in this example we can identify $[-1]$ with the class
$x_1$, and ${\mathcal G}_{F(\sqrt{-1})}=({\mathbb Z}/4)^{n-1}$. Clearly
multiplication by $[-1]$ is injective and so the conclusions of
\ref{7.14} all apply. Indeed if $\zeta_1,\dots ,\zeta_{n-1}$ denote the
regular sequence obtained from the Frattini subgroup, then a basis
for $H^*(S(n))$ as a module over ${\mathbb F}_2[x_1,\zeta_1,\dots ,\zeta_{n-1}]$
is given by the monomials $x_{i_1}x_{i_2}\dots x_{i_s}$ where
$2\le i_1<i_2<\dots <i_s\le n$. In particular the Poincar\'e series
for $H^*(S(n))$ is given by
\[
p_F(t)=\frac{(1+t)^{n-1}}{(1-t)(1-t^2)^{n-1}}=\frac{1}{(1-t)^n}.
\]

We have calculated the cohomology rings of superpythagorean
fields with finite square class group. It is an easy matter to complete our
results and compute $H^*({\mathcal G}_F)$ for any superpythagorean field $F$.
Let $\{[-1]\}\cup\{[a_i]\mid i\in I\}$ be any basis of $\dot F/\dot F^2$ which
contains $[-1]$. Then a set 
$\{\sigma_{-1}\}\cup \{\sigma_i \mid i\in I\}\subset {\mathcal G}_F$ such that
$\sigma_{-1}(\sqrt{-1})=-\sqrt{-1}$, $\sigma_{-1}(\sqrt{a_i})=\sqrt{a_i}$
for all $i\in I$, and $\sigma_i(\sqrt{a_j})=(-1)^{\delta_{ij}}\sqrt{a_j}$
(where $\delta_{ij}=-1$ if $i=j$, $\delta_{ij}=1$ otherwise) and finally
$\sigma_i(\sqrt{-1})=\sqrt{-1}$ is a minimal set of generators for
${\mathcal G}_F$. For any finite subset $S$ of $I$ consider the subgroup
generated by $\sigma_i$, $i\in S$, and by $\sigma_{-1}$. Each $G_S$ is isomorphic
to $S(n+1)$, where $n=|S|$. Moreover each $G_S$ is also a quotient of
${\mathcal G}_F$, and the $(G_S)$ (as $S$ ranges over finite subsets in $I$)
form a projective system of finite subgroups in ${\mathcal G}_F$ such that
${\mathcal G}_F=\varprojlim G_S$.
Therefore from 
\cite[I.2.8]{Se} we see that 
$H^*({\mathcal G}_F)\cong \varinjlim H^*(G_S)$ and hence the
subring generated by one-dimensional classes will simply be
${\mathcal R}={\mathbb F}_2[x_i\mid i\in I]/M$
where $M$ is the ideal generated by the set $\{x_i(x_1+x_i)\mid i\in I\}$.

\end{section}

\begin{section}{Cohomology of Universal W-Groups}
\label{s:universal}

 From the point of view of finite
group cohomology, the universal groups $W(n)$ are perhaps the most interesting
examples of \textit{W}-groups. 
Given that $W(n)$ will surject onto any other \textit{W}-group
${\mathcal G}_F$ for $|\dot F/\dot F^2|=2^n$, $H^*(W(n))$ must
necessarily be rather complicated.
Given that the $W(n)$ satisfy the 2C condition,
we know by \cite{AK} that their 
cohomology is not detectable on proper subgroups,
and hence the usual methods for computing cohomology of finite groups
will run into difficulties. 

Before discussing general properties of these groups, it seems natural
to describe how the universal groups arise as \textit{W}-groups. A basic example
is given by taking the field $F=\mathbb C (t)$, where the elements 
$\{ [t-c] \mid c\in\mathbb C\}$ form a basis for $\dot F/\dot F^2$. The
well-known theorem of Tsen-Lang (see \cite[pp.~45, 46, 296]{Lam1})
implies that each quaternion algebra $(\frac{a,b}{F})$ with 
$a,b\in\dot F$ splits. From theorem 2.20 in \cite{MS1}, we see that
the \textit{W}-group of $F$ is the universal \textit{W}-group on the set
of generators $\{\sigma_{(t-c)}\mid c \in\mathbb C\}$
One can also construct an infinite algebraic extension 
$F\subset K\subset F_q$ (taking $F=\mathbb C (t)$ again) such that
$|\dot K/\dot K^2|=2^n$ ($n$ any prescribed natural number), and
${\mathcal G}_K\cong W(n)$. For details of this construction
we refer the reader to \cite[pp.~102--103]{GMi}.

We now describe some 
subgroups and quotient groups of $W(n)$. 
 From the description we gave for $W(n)$ in terms
of generators and relations, it is apparent that we can construct
a surjective map $\phi_i\colon W(n+1)\to W(n)$ simply by making $x_i=1$
and leaving the other generators unchanged. This map clearly splits
with the canonical injection and hence we have an embedding
$H^*(W(n))\hookrightarrow H^*(W(n+1))$. In fact these maps can
be described in terms of extensions of the form
$$1\to {\mathbb Z}/4\times ({\mathbb Z}/2)^n\to W(n+1)\to W(n)\to 1,$$
where the ${\mathbb Z}/4$ corresponds to the subgroup generated
by $x_i$ and the ${\mathbb Z}/2$ factors correspond to the
commutators $[x_i,x_j]$, where $i\ne j$.

On the other hand, we can also make $n$ of the generators
in $W(n+1)$ equal to $1$, and hence obtain an extension of the
form 
\[
1\to W(n)\times ({\mathbb Z}/2)^n\to W(n+1)\to {\mathbb Z}/4\to 1.
\]

Recall that we have established in \ref{t:cohen-macaulay} that
$H^*(W(n))$ is free and finitely generated over 
the polynomial subalgebra generated by
$r=n + {\binom{n}{2}}$ 2-dimensional classes $\{\zeta_{ij}
\}_{i\le j}$.
Furthermore, the Poincar\'e series for $H^*(W(n))$ is of the
form
\[
p_n(t)=\frac{q_n(t)}{(1-t^2)^{n + \binom{n}{2}}}
\]
where $q_n(t)$ is a palindromic polynomial of degree
$n + {\binom{n}{2}}$ in ${\mathbb Z} [t]$.
Our goal will now be to study the polynomial $q_n(t)$.

Given $E_n$ an elementary abelian 2-group with basis
$\{x_1,\dots ,x_n\}$, we define an action on 
$X_n=(\Sone)^r$, where as before $r=n+{\binom{n}{2}}$.

\begin{definition}
\label{3.6}
Let $(z_{ij})$
denote the complex coordinates for the space $X_n$, ordered
lexicographically. Then
we define the action of $E_n$ on $X_n$ as follows:
\[
x_l(z_{ij})= 
	\begin{cases}
		-z_{ij}, & \text{if $i=j=l$;} \\
		\bar{z}_{ij}, & \text{if $i=l$, $j\ne l$;}\\
		-\bar{z}_{ij}, & \text{if $j=l$, $i\ne l$;}\cr
                z_{ij}, & \text{otherwise.}
\end{cases}
\]
\end{definition}

Following the notation from
\ref{t:X_F}, if ${\mathcal G}_F=W(n)$, we can take
$X_F=X_n$ with the $E_n$ action described as above.
The action is evidently free (indeed, $W(n)$ satisfies
the 2C property) and so the orbit space
$Y_n=X_n/E_n$
is a compact $r$ dimensional
manifold. We shall denote its fundamental group
by $\Upsilon (n)$.
We are now faced with a very explicit topological problem, namely

\begin{problem} 
\label{3.8} 
Calculate $H^*(X_n/E_n)$.
\end{problem}

To begin we consider $W(2)$; it is a group of order 32
expressed as a central extension
$$1\to ({\mathbb Z}/2)^3\to W(2)\to ({\mathbb Z}/2)^2\to 1.$$
Its cohomology has been computed by Rusin \cite{R} but we will
redo this calculation from our point of view. In this
case the group $\Upsilon(2)$ is the fundamental group of
a closed, non-orientable Seifert manifold, i.e.\ a
circle fibering over a 2-dimensional surface (the torus).
Hence it fits into an extension of the form
\[
1\to{\mathbb Z}\to\Upsilon (2)\to{\mathbb Z}\oplus{\mathbb Z}\to 1
\]
with trivial twisting over ${\mathbb F}_2$. This yields a long
exact sequence in mod 2 cohomology,
\begin{multline*}
0\to H^1({\mathbb Z}\oplus{\mathbb Z})\cong H^1(\Upsilon(2))\to \\
\to
H^1({\mathbb Z}) 
\to H^2({\mathbb Z}\oplus{\mathbb Z})\to H^2(\Upsilon(2))\to
H^1({\mathbb Z}\oplus{\mathbb Z},H^1({\mathbb Z}))\to 0
\end{multline*}
as well as the evident isomorphism $H^3(\Upsilon(2))
\cong H^2({\mathbb Z}\oplus{\mathbb Z}, H^1({\mathbb Z}))\cong{\mathbb F}_2$.
Hence we conclude that $q_2(t)=1+2t+2t^2+t^3$.
In terms of ring structure it is easy to verify (using mod 2
Poincar\'e Duality) that
\[
H^*(\Upsilon(2))\cong \Lambda (x_1,y_1,u_2,v_2)
\big\slash \langle R\rangle
\]
where $x_1,y_1,u_2,v_2$ are exterior classes with degree 
equal to their subscript and $R$ is the set of relations
\[
x_1y_1=u_2v_2=x_1u_2=y_1v_2=0, \quad x_1v_2=y_1u_2.
\]
Hence the Poincar\'e series for $W(2)$ is 
\[
p_2(t)=\frac{1+2t+2t^2+t^3}{(1-t^2)^3}.
\]

Now we consider the case $n=3$. Here we have $X_3=(\Sone)^6$
with a free action of $E_3$. If we label these complex
coordinates using $z_{ij}$, then we can project
(equivariantly)
\[
(z_{11}, z_{22}, z_{33}, z_{12}, z_{13}, z_{23})
\mapsto (z_{11}, z_{22}, z_{33}, z_{12}, z_{13})
\]
where $Z=(\Sone)^5$ still has a free $E_3$ action
due to the presence of the 3 free coordinates
$z_{ii}$. Hence we obtain a circle bundle $X_3/E_3\to Z/E_3$.
Similarly, using the projections
\[
(z_{11}, z_{22}, z_{33}, z_{12}, z_{13})
\mapsto (z_{11}, z_{22}, z_{33}, z_{12})
\]
and
\[
(z_{11}, z_{22}, z_{33}, z_{12})\mapsto
(z_{11}, z_{22}, z_{33})
\]
we obtain circle bundles
\[
Z/E_3\to W/E_3 \quad \hbox{and} \quad W/E_3\to U/E_3
\]
where $W\cong (\Sone)^4$ and $U\cong (\Sone)^3$, all
of which have a free $E_3$-action.

For our cohomology computation we start at the bottom,
first we show

\begin{proposition}
\label{4.1}
\[
H^*(W/E_3)\cong\Lambda(u_1,v_1,w_1,b_2,c_2)/ \langle R \rangle
\]
where $R$ is the set of relations
\[
u_1v_1=u_1b_2=v_1c_2=0,\quad u_1c_2=v_1b_2
\]
and the Poincar\'e
series for $W/E_3$ is
$1+3t+4t^2+3t^3+t^4$.
\end{proposition}

\begin{proof} 
First we recall that
$U/E_3\cong (\Sone)^3$. Next we use the Gysin sequence for
the circle bundle $W/E_3\to U/E_3$. Let $e$ represent the
generator for the cohomology of the fiber, and $u_1,v_1,w_1$
generators for the cohomology of the base. The key differential
is given by $d_2(e)=u_1v_1$; the following is an explicit
list of
non-zero classes in $E_3^{*,*}$:
\[
\{ u_1,v_1,w_1,eu_1, ev_1, u_1w_1,v_1w_1, eu_1v_1,
eu_1w_1,ev_1w_1,eu_1v_1w_1\},
\]
from which the calculation above readily follows.
\end{proof}

Our strategy will now be to iterate this method.
We consider the circle bundle $Z/E_3\to W/E_3$;
the Gysin sequence simplifies to yield a sequence
\[
0\to{\mathbb F}_2^3\to H^2(Z/E_3)\to {\mathbb F}_2^3\xrightarrow{d_2}
{\mathbb F}_2^3\to H^3(Z/E_3) 
\to{\mathbb F}_2^3\to 0
\]
and the isomorphisms $H^1(Z/E_3)\cong{\mathbb F}_2^3$,
$H^5(Z/E_3)\cong{\mathbb F}_2$. If $e^{\prime}$ generates
the cohomology of the fiber, then
$d_2(e^{\prime})=u_1w_1$, and hence
$d_2(u_1e^{\prime})=u_1^2w_1=0$, 
$d_2(v_1e^{\prime})=v_1u_1w_1=0$,
$d_2(w_1e^{\prime})=u_1w_1^2=0$
and so $d_2\equiv 0$ in the sequence above. From this
we deduce that $H^*(Z/E_3)$ has Poincar\'e series
equal to
\[
1+3t+6t^2+6t^3+3t^4+t^5
\]
and is in fact generated by 3 one-dimensional classes,
5 two-dimensional classes and one three-dimensional 
class.

The final extension to $X/E_3$ can be done similarly,
using the circle bundle $X/E_3\to Z/E_3$. As before,
the Gysin sequence 
is completely determined by a single differential,
which 
is also identically zero. We obtain

\begin{proposition}
\label{4.2}
The Poincar\'e series for $H^*(X_3/E_3)$ is
\[
q_3(t) = 1 + 3t + 8t^2 + 12t^3 + 8t^4 + 3t^5 + t^6
\]
and hence the Poincar\'e series for $H^*(W(3))$ is
\[
p_3(t)=\frac{1 + 3t + 8t^2 + 12t^3 + 8t^4 +3t^5 + t^6}{(1 - t^2)^6}.
\]
\end{proposition}

Unfortunately this process of iterating Gysin sequences
becomes unmanageable after a few steps. This is related
to interesting questions in differential homological
algebra and homotopy theory which we will discuss further
on.

We will now introduce another torsion-free group
to study the cohomology of \textit{W}-groups. Let $L(n)$ denote
the universal central extension

\begin{equation}
1\to {\mathbb Z}^{\binom{n}{2}}\to L(n)\to {\mathbb Z}^n\to 1
\end{equation}
which is often called the free 2-step nilpotent group on 
$n$ generators. Indeed the group $L(n)$ is generated by
elements $\sigma_1,\sigma_2,\dots ,\sigma_n$ of infinite order
with the relations imposed by requiring that the commutators
$[\sigma_i,\sigma_j]$ for $i<j$ be central (hence the extension
above).

Recall from our description of $W(n)$ that it
is generated by elements $x_1,\dots ,x_n$ of order $4$ such
that their squares and commutators are central (and of order 2).
Hence if we abelianize $W(n)$ we obtain $({\mathbb Z}/4)^n$, and so
a central extension

$$1\to ({\mathbb Z}/2)^{\binom{n}{2}}\to W(n)\to ({\mathbb Z}/4)^n\to 1$$
where the generators of the kernel correspond to the commutators.
Reducing mod $4$ we obtain a map $L(n)\to ({\mathbb Z}/4)^n$ which
factors through $W(n)$ and hence we obtain an extension
\begin{equation}
\label{5.2}
1\to L_4(n)\to L(n)\to W(n)\to 1
\end{equation}
where $L_4(n)$ appears as some sort of ``congruence subgroup''.
Note that it will be generated by the fourth powers
$u_i=\sigma_i^4$ and the squares of the 
commutators $w_{ij}=[\sigma_i,\sigma_j]^2$. 

We begin our cohomological analysis by computing the mod 2
cohomology ring of $L_4(n)$.

\begin{theorem}
\label{5.3}
The mod 2 cohomology of $L_4(n)$ is an exterior algebra
on $r=n + {\binom{n}{2}}$ one dimensional generators.
\end{theorem}
\begin{proof}
 From the above it is clear that we can also realize $L_4(n)$
as a central extension of the subgroup generated by commutators,
\[
1\to {\mathbb Z}^{\binom{n}{2}}\to L_4(n)\to {\mathbb Z}^n\to 1.
\]
To compute the mod 2 Lyndon-Hochschild-Serre spectral sequence
for this extension we need to determine the commutators 
$[u_i,u_j]$ for $i<j$. A straightforward calculation shows that
$\sigma_i^4\sigma_j^4\sigma_i^{-4}\sigma_j^{-4}
=[\sigma_i,\sigma_j]^{16}$. Now if we consider the $E_2$ term
of the aforementioned spectral sequence, it is of the form
$\Lambda (u_1,u_2,\ldots ,u_n)\otimes\Lambda (v_1,\ldots, 
v_{\binom{n}{2}})$
where by abuse of notation
we have identified generators from the base with their duals.
 From the extension data computed above, we see that the $v_1,\ldots
,v_{\binom{n}{2}}$ are permanent cocycles, as they transgress to elements
divisible by 8 and hence 0 mod 2. The spectral sequence 
collapses at $E_2$, there are evidently no extension problems
and the result follows.
\end{proof}

We will now use this to analyze the cohomology of $L(n)$.

\begin{theorem}
\label{5.4}
In the mod 2 LHS spectral sequence associated to 
\[
1\to L_4(n)\to L(n)\to W(n)\to 1
\]
the one dimensional generators in the cohomology of $L_4(n)$
transgress to a regular sequence in $H^2(W(n),{\mathbb F}_2)$,
$E_3=E_{\infty}$, and in particular if $I$ is the ideal 
generated by these transgressions, then  
\[
H^*(W(n),{\mathbb F}_2)/I
\cong H^*(L(n),{\mathbb F}_2).
\]
\end{theorem}

\begin{proof}
Let $L_2(n)$ denote the subgroup of $L(n)$ generated by
the commutators $w_{ij}$, $i<j$ and the squares $\sigma_i^2$;
then this group fits into a commutative diagram of
extensions:
\[
\xymatrix{
& & 1 & 1 \\
& & ({\mathbb Z}/2)^n \ar@{=}[r] \ar[u] &  ({\mathbb Z}/2)^n \ar[u] \\
1 \ar[r] & L_4(n) \ar[r]  & L(n) \ar[u] \ar[r]	& W(n) \ar[r] \ar[u] &	1 \\
1 \ar[r] & L_4(n) \ar[r] \ar@{=}[u] & L_2(n) \ar[r] \ar[u] & ({\mathbb Z}/2)^r \ar[r] \ar[u] & 1 \\  
& & 1 \ar[u] & 1 \ar[u] \\
}
\]
Note that $H=({\mathbb Z}/2)^r\subset W(n)$ is the maximal central
elementary abelian subgroup. For reasons analogous to the
arguments in \ref{5.3}, $H^*(L_2(n),{\mathbb F}_2)$ is also an exterior
algebra on $r$ one dimensional generators. Consider the
mod 2 LHS spectral sequence for the bottom row. Let
$I^{\prime}$ denote the ideal generated by the transgressions.
Then evidently $H^*(({\mathbb Z}/2)^r)/I^{\prime}\cong H^*(L_2(n))$,
and so we conclude that the ideal is precisely the ideal generated
by the squares of the one dimensional generators. Comparing
spectral sequences, we see that $I$ is generated by trangressions
$\zeta_{ij}$, $i<j$ such that they restrict to squares in
$H^2(H,{\mathbb F}_2)$. As we saw in section~\ref{s:cohomology} these
form a regular sequence 
of maximal length in the cohomology of $W(n)$. Indeed, they
generate a polynomial
subalgebra over which it is free and finitely generated. From
this it follows that $E_3=E_{\infty}$ and hence that
$H^*(L(n),{\mathbb F}_2)\cong H^*(W(n))/(\zeta_{ij})$.
\end{proof}

\begin{remark} From the above we see that the mod 2 cohomology
of $L(n)$ plays the same role as the cohomology of $\Upsilon (n)$.
An important difference is that the latter group is virtually
abelian, whereas $L(n)$ is not. The part played by the associated
circle classes in cohomology is played by the cohomology classes
in $H^1(L_4(n),{\mathbb F}_2)$. In other words we use a ``mod 2
cohomology torus'' instead. Note that the integral cohomology of $L_4(n)$ can be
quite complicated.
\end{remark}

The interesting feature of $L(n)$ is that it has substantial
torsion-free cohomology. Furthermore, there are reasons to believe that
all of it may be 2-torsion-free. In contrast, $H^*(\Upsilon (n),\mathbb Q)$
is fairly small, even though it has the same mod 2 cohomology as
$L(n)$. Later we shall provide a complete calculation for the rational
cohomology $H^*(L(n),\mathbb Q)$. From our constructions it is clear
that these rational classes come from integral classes which in turn
will appear as mod 2 classes in the cohomology of $W(n)$. This seems
like a rather novel approach---using rational cohomology to produce 
mod 2 cohomology for a finite group.

We will now describe a well-known topological method for approaching
the cohomology of $L(n)$. This is based on \cite{GM} and \cite{L}, although
these methods were known even before these references appeared.
To begin we observe that the central extension defining $L(n)$
gives rise to
a classifying map $\phi\colon (\Sone)^n\longrightarrow 
(\mathbb C P^{\infty})^{\binom{n}{2}}$. If we consider the pullback of
the universal bundle over $(\mathbb C P^{\infty})^{\binom{n}{2}}$, we obtain
(up to homotopy) a diagram of fibrations
\[
\label{5.5}
\begin{CD}
BL(n) @>>> U \\
@VVV @VVV \\
(\Sone)^n @>{\phi}>> ({\mathbb C}P^{\infty})^{\binom{n}{2}} \\
\end{CD}
\]
where $BL(n)$ is the classifying space for $L(n)$ (and hence its 
cohomology is the group cohomology of $L(n)$) and the right hand
column is the universal bundle (and so $U$ is contractible).
Note that the map $\phi$ is determined by the extension data, which
in this case is the fact that the commutators are central. Note
that if $b_{ij}$ for $i<j$ form an integral basis for
$H^2((\mathbb C P^{\infty})^{\binom{n}{2}},{\mathbb Z})$, then we can 
assume $\phi^*(b_{ij})=x_ix_j$, where $x_1,\dots ,x_n$ is an
exterior basis of one dimensional elements for the cohomology
of $(\Sone)^n$.

Associated to any such pull-back diagram of spaces 
we have an \textsl{Eilenberg-Moore spectral sequence}
converging to the cohomology of the upper
left hand corner 
(in this case the cohomology of $L(n)$). In our situation
the $E_2$ term can be
computed homologically as
\[
\label{5.6}
E_2^{*,*}=\Tor_{{\mathbb Z}[b_{ij}]}(\Lambda
(x_1,\dots ,x_n), {\mathbb Z}),
\]
where we have written $\Lambda(x_1,\dots ,x_n)$ for $H^*((\Sone)^n,
{\mathbb Z})$.  
\begin{remark}
This identification means that our convention for exterior algebras
over ${\mathbb Z}$ requires the relation $e\wedge e = 0$ and not just
$e \wedge e = - e \wedge e$.  In particular, $\Lambda(x_1,\dots ,x_n)$
is a free ${\mathbb Z}$-module. 
\end{remark}

We can obtain a great deal of information about this $E_2$ term using
methods from commutative algebra. The map $\phi$ induces
an algebra map
\[
\phi^*\colon {\mathbb Z}[b_{ij}]\longrightarrow \Lambda (x_1,x_2,\dots ,x_n)
\]
with $\phi^*(b_{ij})=x_ix_j$. To compute the $\Tor$ term in the
spectral sequence, we make use of a particular free ${\mathbb Z}[b_{ij}]$
resolution for ${\mathbb Z}$, namely the Koszul complex (see \cite{L} for
complete details) which is of the form
$\Lambda (u_{ij})\otimes {\mathbb Z}[b_{ij}]$. 
The $E_2$ term above can then be computed in the usual way, i.e.\
setting 
\begin{equation}
\label{e:K-definition}
K^{*,*} = \Lambda (u_{ij})\otimes {\mathbb Z}[b_{ij}]
\otimes_{{\mathbb Z}[b_{ij}]}
\Lambda(x_1,\dots ,x_n),  
\end{equation}
we have $E_2^{*,*} = H(K^{*,*})$.  In fact, more is true: the EMSS
starts out with an $E_1$-term which can be indentified with our
complex $K^{*,*}$. 

Let us now record a few facts about the functoriality of this complex
which will be used later.  First, we have been working with integer
coefficients in cohomology; if we liked we could tensor $K^{*,*}$ with
${\mathbb Q}$ or ${\mathbb F}_2$ to get the $E_1$-term for the
rational or mod 2 version of the EMSS.  Second, note that $K^{*,*}$
can be written in coordinate-free fashion as $\Lambda^*\Lambda^2 V
\otimes \Lambda^*V$, where $V$ is a free ${\mathbb Z}$-module of rank
$n$.  In the bigrading associated with this description (which is
\emph{not} the same as the usual bigrading for the EMSS), the
differential has bidegree $(2,-1)$ and a simple coordinate-free
description: if $\theta \in \Lambda^r V$, 
and $\omega_1, \dots, \omega_m \in \Lambda^2 V$, then
\begin{equation}
\label{e:diff}
d_1 (\theta \otimes \omega_1 \wedge \dots \wedge \omega_m) = 
\sum_{i=1}^{m} (-1)^{i+1} (\theta \wedge \omega_i) \otimes (\omega_1
\wedge \dots \wedge \widehat{\omega_i} \wedge \dots \wedge \omega_m),
\end{equation}
where the sum on the right hand side is an element of $\Lambda^{r+2} V
\otimes \Lambda^{m-1} \Lambda^2 V$. These coordinate-free descriptions
show that $K$ is functorial in $V$; in particular $K \otimes {\mathbb
Q}$ and its homology are representations of $GL(V \otimes {\mathbb
Q})$.  

At this stage we have two problems to address. First we have
the purely algebraic problem of computing the homology of
the Koszul complex. As we shall see shortly this has been done
for rational coefficients, but for coefficients in $\mathbb F_2$ it is 
an open problem. Secondly, to make this an effective method
of computation we need to determine the higher differentials
or whether in fact $E_2=E_{\infty}$. If the combinatorial
determination of the homology of the Koszul complex and the
collapse at $E_2$ can \emph{both} be established then we will
have a complete computation. We should also note that the
$E_3$ term of the Lyndon-Hochschild-Serre spectral sequence
for the central extension defining $L(n)$
will also be isomorphic to the cohomology
of the Koszul complex. Hence its collapse at $E_3$ is equivalent
to the collapse of the EMSS at $E_2$. There are however some 
advantages in dealing with higher differentials in the EMSS
related to methods from homological perturbation theory
(see \cite{GM}).

Over the rationals there is a very nice collapse theorem.

\begin{proposition}
\label{5.9}
The Eilenberg-Moore spectral sequence 
for $R=\mathbb Q$ collapses at $E_2=E_{\infty}$ without extension
problems, and hence we have
\[
H^*(L(n),\mathbb Q)\cong H(K \otimes {\mathbb Q}).
\]
\end{proposition}
\begin{proof} This is proposition 4.3.1 in \cite{L}.
Without elaborating too much we shall only say that this
works because we can use ``rational de Rham complexes'' which
are commutative, from which it follows that the higher
differentials are zero.
\end{proof}

Note that the isomorphism above is an isomorphism of algebras,
i.e.\ we can endow the homology of the Koszul complex with a
natural product induced from that on the original complex.

Combining \ref{5.4} with \ref{5.9} we can construct a large number of 
mod 2 cohomology classes for $W(n)$. Let $T$
denote the ideal of torsion classes, then we have an inclusion
\[
H^*(L(n),{\mathbb Z})/T \otimes{\mathbb F}_2
\hookrightarrow H^*(W(n))/(\zeta_{ij}).
\]
Note that in particular we can determine the dimension of this
subspace as the dimension of the \emph{rational} cohomology of
$L(n)$ or of the corresponding rational Koszul complex.

It is interesting to link the ``algebraic 2-torsion'' in the
Koszul complex with 2-torsion in the cohomology of $L(n)$.
Both of these are open questions as far as we know. They
can be easily related via 
\begin{proposition}
\label{5.10}
The cohomology $H^*(L(n),{\mathbb Z})$ has no 2-torsion and the
associated mod 2 Eilenberg-Moore spectral sequence collapses
at $E_2$ if and only if the homology of the Koszul complex
$K^{*,*}$ has no 2-torsion.
\end{proposition}
\begin{proof} Suppose that $H^*(L(n),{\mathbb Z})$ has no
2-torsion and that the mod 2 EMSS collapses at $E_2$. This means
that this $E_2$ term accounts for all the mod 2 cohomology, which
by our hypothesis must have the same dimension as the rational
cohomology. In terms of the Koszul complex this simply means
that the homology of the rational version must have the same
total dimension as the mod 2 complex. Hence we know that the
homology of the integral complex must be 2-torsion free. Conversely
if this condition holds, the $E_2$ terms of both the rational and
mod 2 EMSS for $L(n)$ must have the same dimensions, and as the
rational one always collapses, the same must hold for the mod 2
one and in addition there cannot be any 2-torsion in $H^*(L(n),{\mathbb Z})$,
as this would produce unaccountable mod 2 cohomology classes.
\end{proof}

An explicit calculation due to Lambe \cite{L} shows that for $n\le 5$
the conditions in \ref{5.10} hold. Hence we have complete calculations
for the cohomology of the corresponding \textit{W}-groups. In other words
if as before $I\subset H^*(W(n),{\mathbb F}_2)$ is the ideal generated
by the regular sequence $\zeta_{ij}$, then we have an extension 
of algebras
\[
0\longrightarrow I\longrightarrow H^*(W(n),{\mathbb F}_2)
\longrightarrow H(K\otimes {\mathbb F}_2)
\longrightarrow 0.
\]
Using the computations in \cite{L} we can now record the Poincar\'e
series for $W(4)$ and $W(5)$, following the format previously established.
We have that $p_4(t)=q_4(t)/(1-t^2)^{10}$, where $q_4(t)$ is the
polynomial
\[
{1+4t+20t^2+56t^3+84t^4+90t^5+84t^6+56t^7+20t^8+4t^9+t^{10}}.
\]
For $p_5(t)$ the denominator will be $(1-t^2)^{15}$ and the numerator
will be $q_5(t)$, the polynomial
\begin{gather*}
1+5t+40t^2+176t^3+440t^4+835t^5+1423t^6+1980t^7 \\
+1980t^8+1423t^9
+835t^{10}+440t^{11}+176t^{12}+40t^{13}+5t^{14}+t^{15}.
\end{gather*}

We will now describe the homology of the rational
Koszul complex defined previously. The computation of $H(K \otimes
{\mathbb Q})$ has appeared in many guises (see \cite{JW}, also
\cite{Bouc}, \cite{Ka}, and \cite{Si}), so we do not repeat it here.
The connections between the various contexts in which this computation
occurs are not completely transparent, so we point out that the
calculation is equivalent to computing the Lie Algebra cohomology of   
the graded Lie algebra associated to the group $L(n)$.  To state the
result, let us take $V$ to be an $n$-dimensional rational vector space,
identify $K \otimes {\mathbb Q}$ with $\Lambda^*\Lambda^2 V
\otimes \Lambda^*V$, and use the associated bigrading. As noted
earlier, this complex and its homology are representations of $GL(V)$.
Irreducible representations of this group are parametrized by Young
diagrams, which are a way of representing partitions of natural
numbers; the reader unfamiliar with this theory may wish to consult
\cite{FH}.

\begin{theorem}
\label{6.1}
The homology of $K\otimes\mathbb Q$ in the $(p,q)$ position is the sum of all
representations corresponding to symmetric $p+2q$-box $p$-hook
diagrams.
\end{theorem}

The corresponding rational ranks are easily computable by a standard
dimension formula \cite{McD}. We state this formula for the
convenience of the reader. Let $\lambda = (\lambda_1, \ldots,
\lambda_n)$ be a partition with trailing zeroes added to fill out its
length to $n$ if necessary. Recall that by
convention, we take $\lambda_1 \geq \lambda_2 \geq
\cdots \geq \lambda_n$. 
Associated to $\lambda$ we have a diagram $Y_{\lambda}$
with $\lambda_1$ boxes
in the first column, $\lambda_2$ in the second, etc.\ where the
columns all begin on the same horizontal line. Hence our partition
gives a diagram with $n$ columns of boxes, of length $\lambda_1,\dots,
\lambda_n$ respectively (some at the end may be zero). In total we
will have $\lambda_1 +\dots + \lambda_n$ boxes
in $Y_{\lambda}$. We label these
boxes by pairs $(i,j)$ corresponding to rows and columns. 
Note that a diagram will be \textsl{symmetric} if it is invariant 
under the \textsl{transposition} of exchanging rows and columns. The number
of ``hooks'' will be equal to the number of boxes in the diagram which
are bisected by the diagonal; to each box we associate a
\textsl{hook} by taking its union with all boxes below it and all
boxes to the right of it.  By definition, the
\textsl{hooklength} of the box $(i,j)$ is 
$h(i,j)$, the number of boxes
below $(i,j)$ plus  the number of boxes to the right of $(i,j)$
plus $1$. 
Let $V$ be a vector space of dimension
$n$ and let $S_{\lambda}V$ be the representation of $GL(V)$ corresponding
to $\lambda$, then we have:
\[
\label{6.2}
\dim S_{\lambda}V = \prod_{(i,j)\in Y_{\lambda}}
\frac{n+j-i}{h(i,j)}
\]
We can use this formula to find a lower bound on the coefficients
of the Poincar\'e series $p_n(t)$. Indeed we have
\begin{theorem}
\label{t:combinatorics}
Let $p_n(t)=q_n(t)/(1-t^2)^r$ denote the Poincar\'e series
for $H^*(W(n))$, where $r=\binom{n+1}{2}$. If 
$q_n(t)=1+a_1t+a_2t^2+\dots + a_rt^r$, then
\[
a_k\ge \sum_{p+q=k} \sum_{Y_{\lambda}} \prod_{(i,j)\in Y_{\lambda}}
\frac{n+j-i}{h(i,j)}
\]
where $Y_{\lambda}$ ranges over all symmetric $p+2q$-box, $p$-hook
Young diagrams, and $h(i,j)$ denotes the hooklength of the box
(i,j).
\end{theorem}

As we have noted, the rational cohomology $H^{p,q}(K \otimes {\mathbb
Q})$ obtained from Young diagrams gives us mod 2 cohomology in
$H^{p,q}(K \otimes {\mathbb F}_2)$.  If we show that the
\emph{integral} cohomology $H^{p,q}(K)$ is 2-torsion free, we will be
able to conclude that there is no
mod 2 cohomology unaccounted for by our constructions, or in other
words, that the inequality of theorem~\ref{t:combinatorics} is actually an
equality.  We will demonstrate by ad hoc methods that this
is the case for certain small values of $p$ and $q$: 
\begin{lemma}
\label{l:koszul-torsion}
$H^{p,q}(K)$ is 2-torsion free if $p+q \leq 4$ and $q \leq 3$. 
\end{lemma}
It then follows from  the universal coefficient theorem that we have:
\begin{theorem}
\label{6.3}
For all integers $n\ge 1$, we have 
\[
a_2=\frac{n(n+1)(n-1)}{3} \quad \text{and} \quad
a_3=\frac{n(n^2-1)(3n-4)(n+3)}{60}. 
\]
\end{theorem}
These results could certainly be extended, but as
we have no method that would give all the coefficients $a_k$, we have
elected to stop here. 

This section has been rather technical, so for the weary reader who
would like to see some easily-understandable consequences of the
results we offer the following facts. 
 From the very definition of the groups $W(n)$ we know that their
one dimensional cohomology classes multiply trivially, i.e.\ the
subring in ${\mathcal R}\subset H^*(W(n))$ generated by $H^1(W(n))$
is an $n$ dimensional space with no products. Hence 
there are precisely 
$n(n+1)(2n+1)/6$
2-dimensional
generators in $H^*(W(n))$.
More generally, if $n$ is sufficiently large,
the $k$-th coefficient $a_k$ of
$q_n(t)$
will satisfy the inequality
\[
a_k\ge \frac{(n+k-1)(n+k-2)\cdots (n-k+1)}{(2k-1)[(k-1)!]^2},
\]
which can be derived from theorem~\ref{t:combinatorics}.

We finish this section with the proof of our lemma about the 2-torsion
freeness of $H^{p,q}(K)$. 
\begin{proof}[Proof of lemma~\ref{l:koszul-torsion}]
Recall that the groups $K^{p,q} = \Lambda^p V \otimes \Lambda^q
\Lambda^2 V$ are free abelian, and that the differential
(\ref{e:diff}) has bidegree $(2,-1)$.  If $p$ is $0$ or $1$, then
$H^{p,q}$ is a kernel and therefore also free abelian, in particular
$2$-torsion free.  It is easy to show that $H^{p,q}$ is zero if $q =
0$ and $p > 1$, so we need only consider the cases $(p,q) = (2,2)$,
$(2,1)$, and $(3,1)$. 
\begin{proof}[$(p,q) = (2,1)$]
First note that for any ${\mathbb Z}$-module $W$, there is an exact
sequence $\Lambda^2 W \to W \otimes W \to S^2 W$, where the first map
sends $u \wedge v$ to $u \otimes v - v \otimes u$. Applying this to $W
= \Lambda^2 V$, we see that the cokernel of $\Lambda^2 \Lambda^2 V \to
\Lambda^2 V \otimes \Lambda^2 V$ is $S^2 \Lambda^2 V$, so $H^{2,1}$ is
a submodule of the free module $S^2 \Lambda^2 V$, and therefore
2-torsion free. 
\renewcommand{\qedsymbol}{}
\end{proof}
\begin{proof}[$(p,q) = (3,1)$]
We show that the sequence 
\[
V \otimes \Lambda^2 \Lambda^2 V \to^\phi \Lambda^3 V \otimes \Lambda^2 V
\to^\mu \Lambda^5 V
\]
is exact at $\Lambda^3 V \otimes \Lambda^2 V$, so $H^{3,1}$ is zero. 
To see this, let $\{e_i\}$ be a basis for $V$, note that the
kernel of $\mu$ is generated by elements of the form $e_i \wedge e_j
\wedge e_k \otimes e_l \wedge e_m - e_i \wedge e_l
\wedge e_m \otimes e_j \wedge e_k$, and that all of these elements are
in the image of $\phi$.
\renewcommand{\qedsymbol}{}
\end{proof}
\begin{proof}[$(p,q) = (2,2)$]
Notice that the natural map $\mu\colon \Lambda^2 V \otimes \Lambda^2
\Lambda^2 V \to \Lambda^3 \Lambda^2 V $ given by $\omega_1 \otimes
\omega_2 \wedge \omega_3 \mapsto \omega_1 \wedge \omega_2 \wedge
\omega_3$ has $\mu d_1 = 3$ (see equation~\ref{e:diff}). For purposes
of 2-torsion, this means that the differential $\dif_1\colon \Lambda^3
\Lambda^2 V \to \Lambda^2 V \otimes \Lambda^2 \Lambda^2 V$ is split,
and that there is therefore no 2-torsion in $H^{2,2}$.  
\renewcommand{\qedsymbol}{}
\end{proof}
This completes the proof of lemma~\ref{l:koszul-torsion}. 
\end{proof}

\end{section}

\begin{section}{Final Remarks}
\label{s:final remarks}
In the preceding sections we have described the basic cohomological
structure of \textit{W}-groups by making use of certain topological
models. However one could equally well attempt to compute 
$H^*({\mathcal G}_F)$ directly from the central extension
\[
1\to\Phi ({\mathcal G}_F)\to {\mathcal G}_F\to E\to 1.
\]
The most meaningful situation occurs when $|\dot F/\dot F^2|=2^n$,
in which case $E=E_n\cong ({\mathbb Z}/2)^n$
and $\Phi ({\mathcal G}_F)\cong ({\mathbb Z}/2)^r$. There is an
Eilenberg-Moore spectral sequence associated to this extension, with
$E_2$-term given by the bigraded algebra
\[
\hbox{Tor}_{H^*(K(\Phi ({\mathcal G}_F),2))}(H^*(E_n), {\mathbb F}_2).
\]
Now if $V$ is an elementary abelian $2$-group, $H^*(K(V,2))$
is a polynomial algebra on countably many generators.
By using the explicit form of the $k$-invariants as described
in \ref{t:central} this algebra can be simplified as follows. Let
$\kappa_1,\dots ,\kappa_r\in H^2(E_n)$ denote the $k$-invariants of
the extension. They can be used to define a map of polynomial algebras:
\[
{\mathbb F}_2 [b_1,\dots ,b_r] \to H^*(E_n)
\]
where the $b_1,\dots ,b_r$ have degree $2$ and $b_i\mapsto \kappa_i$.
Then we can express the $E_2$ term above as an extension
\[
0\to (\zeta_1,\dots ,\zeta_r)\to
E_2^{*,*}\to \hbox{Tor}_{{\mathbb F}_2[b_1,\dots ,b_r]}({\mathbb F}_2
[\dot F/\dot F^2], {\mathbb F}_2)\to 0
\]
where the polynomial classes
$\zeta_1,\dots ,\zeta_r$ are permanent cocyles in bidegree
$(-1, 3)$, and can be chosen to represent
the regular
sequence we have already obtained (see \ref{t:cohen-macaulay}).
We now make a conjecture which has been verified for all examples
we know
\begin{conjecture}
Up to filtration, we have an isomorphism of algebras
\[
H^*({\mathcal G}_F)/(\zeta_1,\dots ,\zeta_r)\cong
\hbox{Tor}_{{\mathbb F}_2[b_1,\dots ,b_r]}({\mathbb F}_2[\dot F/\dot F^2],
{\mathbb F}_2)
\]
\end{conjecture}

The validity of this conjecture remains an interesting open
question. It is equivalent to the collapse at $E_3$
for the Lyndon-Hochschild
Serre spectral associated to the extension above.
The advantage of the EMSS is the fact that in this kind of
situation models with
explicit differentials have been developed (see \cite{GM})
and hence substantial insight can be obtained; however a definitive proof 
(or counterexample) would seem to require additional ideas.
It is interesting to note however, that the Galois cohomology
occurs as an edge in the spectral sequence above, and the
Tor algebra we describe seems to be the most natural
extension of this to a global computation for \textit{W}-groups.
As we have mentioned previously, 
collapse on the edge is implied by the Milnor Conjecture.
It would seem reasonable to expect that the intrinsic
field theory input which determines ${\mathcal G}_F$
will play an important role here. The
higher differentials can be determined using cup-$1$ products;
one can only speculate that a sensible approach would be
to get a hold on them by using the field theory context in which
\textit{W}-groups are defined.

\end{section}


\begin{thebibliography}{99}

\bibitem{AK} Adem, A. and Karagueuzian, D. \emph{Essential
Cohomology of Finite Groups}, Comm. Math. Helv.
\textbf{72} (1997), pp.~101--109.

\bibitem{AKM} Adem, A., Karagueuzian, D. and Minac, J. \emph{Topological
Models and the Cohomology of Galois Groups}, Comptes Rendus Acad. Sci. Paris,
t. 326, Serie I (1998), pp.~919--324.

\bibitem{AM} Adem, A. and R. J. Milgram, R. J. \emph{Cohomology of
Finite Groups} (Springer-Verlag Grundlehren 309), Springer-Verlag,
Berlin, 1994.

\bibitem{AP} Allday, C. and Puppe, V. \emph{Cohomological Methods in
Transformation Groups}, Cambridge University Press, 1994.

\bibitem{BC} Benson, D. and  Carlson, J. \emph{Projective Resolutions
and Poincar\'e Duality Complexes}, Transactions of the A.M.S.
\textbf{342} (1994), pp.~447--488.

\bibitem{Bouc} Bouc, S. \emph{Homologie de certains ensembles de $2$-sous-groupes des groupes sym\'etriques}.
(French) [Homology of certain sets of $2$-subgroups of symmetric groups] J. Algebra \textbf{150} (1992), no. 1,
pp.~158--186. 

\bibitem{Bredon} Bredon, G. \emph{Introduction to compact
transformation groups.} Pure and Applied Mathematics, Vol. \textbf{46}.
Academic Press, New York-London, 1972.  

\bibitem{Bro} Br\"ocker, L. \emph{Uber die Anzahl der Anordnungen
Eines kommutativen K\"orpers}, Arch. Math. \textbf{29} (1977), pp.~149--163.

\bibitem{BH} Broto, C. and  Henn, H.-W. \emph{Some remarks on central
elementary abelian $p$-subgroups and cohomology of classifying
spaces}, Quart. J. Math. Oxford Ser. (2) \textbf{44} (1993), no. 174, 
pp.~155--163. 

\bibitem{C1} Carlson, J. 
\texttt{http://www.math.uga.edu/\symbol{126}jfc/groups/Groups64}. 

\bibitem{C2} Carlson, J. \emph{Depth and Transfer Maps in the Cohomology
of Groups}, Math. Zeit. \textbf{218} (1995), pp.~461--468.

\bibitem{Ca} Carlsson, G. \emph{On the Non--existence of Free Actions
of Elementary Abelian Groups on Products of Spheres},
Am. J. Math. \textbf{102} (1980), pp.~1147--1157.

\bibitem{Ch} Charlap, L. \emph{Bieberbach Groups and Flat Manifolds},
(Springer-Verlag Universitext) Springer-Verlag, Berlin, 1986.

\bibitem{CS} Craven, T. and Smith, T. \emph{Formally Real Fields from
a Galois Theoretic Perspective}, to appear in Journal of Pure and
Applied Algebra.

\bibitem{D} Duflot, J. \emph{Depth and Equivariant Cohomology},
Comm. Math. Helv. \textbf{56} (1981), pp.~627--637.

\bibitem{EGH} Eisenbud, D., Green, M., and Harris, J.
\emph{Cayley-Bacharach Theorems and Conjectures}, Bull. 
Amer. Math. Soc. (N.S.) \textbf{33} (1996), no. 3, pp.~295--324. 

\bibitem{E} Evens, L. \emph{Cohomology of Groups}, Oxford University
Press, Oxford, 1991.

\bibitem{FH} Fulton, W. and Harris, J. \emph{Representation
Theory, a First Course} (GTM 129), Springer-Verlag, Berlin, 1991.

\bibitem{G} Gao, W. Ph.D. Thesis,
U. Western Ontario, 1996. 

\bibitem{GMi} Gao, W.  and Minac, J. \emph{Milnor's Conjecture and
Galois Theory I}, Fields Institute Communications Vol.~16 (1997),
pp.~95--110. 

\bibitem{GM}  Gugenheim, V.K.A.M.  and May, J.P. \emph{On the Theory and 
Application of Torsion Products},
Memoirs of the A.M.S. \textbf{142} (1974).

\bibitem{H} Hillmann, J. \emph{Flat 4--manifold Groups},
New Zealand J. Math. \textbf{24} (1995), pp.~29--40.

\bibitem{HS} Hilton, P. and Stammbach, U. \emph{An Introduction
to Homological Algebra} (GTM 4), Springer Verlag, Berlin, 1971.

\bibitem{JW} Jozefiak, T. and Weyman, J. \emph{Representation--theoretic
interpretation of a formula of D.E. Littlewood},
Math Proc. Cambridge Phil. Soc. \textbf{103} (1988), pp.~193--196.

\bibitem{Ka} Karagueuzian, D. \emph{Homology of Complexes of Degree One
Graphs}, Ph.D. Thesis, Stanford University,  1994.

\bibitem{Lam1} Lam, T. Y. \emph{The Algebraic Theory of Quadratic
Forms}, Benjamin/Addison--Wesley, Reading, Mass., 1973.

\bibitem{Lam2} Lam, T. Y. \emph{Orderings, Valuations and Quadratic
Forms}, C.B.M.S. Vol.~52, A.M.S. (1983).

\bibitem{L} Lambe, L. \emph{Cohomology of Principal $G$--Bundles over
a Torus when $H^*(BG,R)$ is Polynomial},
Bull. Soc. Math. Belgique \textbf{38} (1986), pp.~247--264.

\bibitem{La} Lang, S. \emph{Algebra} (Third Edition), Addison-Wesley,
Reading, Mass., 1993.

\bibitem{Laz} Lazard, \emph{Groupes Analytiques p-adiques},
Inst. Hautes Etudes Scientifiques, Publ. Math. \textbf{26}, pp.~389--603.

\bibitem{McD} Macdonald, I. G. \emph{Symmetric functions and Hall
polynomials.} Second edition. With contributions by A. Zelevinsky.
Oxford Mathematical Monographs. Oxford Science Publications. The 
Clarendon Press, Oxford University Press, New York, 1995.

\bibitem{Mer} Merkurjev, A. \emph{On the Norm Residue Symbol of
Degree Two}, Doklad. Akad. Nauk. SSR \textbf{261} (1981), pp.~542--547
(English Translation: Soviet Math. Doklady \textbf{24} (1981) pp.~546--551.)

\bibitem{M} Milnor, J. \emph{Algebraic \textit{K}-theory and Quadratic Forms},
Inventiones Math. \textbf{9} (1970), pp.~318--344.

\bibitem{MMS} Mah\'e, L.,  Minac, J.,  and Smith, T. \emph{Additive Structure
of Subgroups of $\dot F/\dot F^2$ and Galois Theory}, preprint 1997.

\bibitem{Mi} Minac, J. \emph{Elementary $2$-Abelian Subgroups of
\textit{W}-groups}, preprint. 

\bibitem{MSm1} Minac J.  and Smith, T. \emph{Decomposition of Witt Rings
and Galois groups},  Can. J. Math. \textbf{47} (1995), pp.~1274--1289.

\bibitem{MSm2} Minac, J. and Smith, T. \emph{\textit{W}-groups of Quadratic
Extensions}, in preparation.

\bibitem{MS1} Minac, J. and Spira, M. \emph{Witt Rings and Galois Groups},
Annals of Mathematics \textbf{144} (1996), pp.~35--60.

\bibitem{MS2} Minac, J. and Spira, M. \emph{Formally Real Fields,
Pythagorean Fields, C-fields and \textit{W}-groups}, Math. Zeit.
\textbf{205} (1990), pp.~519--530.

\bibitem{OVV} Orlov, D., Vishik, A. and Voevodsky V. \emph{Motivic
Cohomology of Pfister Quadrics}, in preparation.

\bibitem{Ribes} Ribes, L. \emph{Introduction to Profinite Groups and
Galois Cohomology}, Queen's Papers in Pure and Applied Mathematics
No. 24, Queen's University, Kingston, Ontario Canada (1970).

\bibitem{R} Rusin, D. \emph{The Cohomology of the Groups of Order 32},
Math. Comp. \textbf{53} No.~187 (1989), pp.~359--385.

\bibitem{Serre-loc} Serre, J.-P. \emph{Alg\`ebre
Locale---Multiplicit\'es}, Springer Lecture Notes in Mathematics,
\textbf{11}, Springer-Verlag, Berlin, 1965. 

\bibitem{Se} Serre, J-P. \emph{Galois Cohomology}, Springer-Verlag,
Berlin, 1997. 

\bibitem{Si} Sigg, S. \emph{The Laplacian and the Homology of Free
Two-step Nilpotent Lie Algebras}, J. Algebra \textbf{185} (1996), pp.~144--161.

\bibitem{V} Voevodsky, V. \emph{The Milnor Conjecture}, preprint 1997.

\bibitem{W} Wolf, J. \emph{Spaces of Constant Curvature} (5th edition),
Publish or Perish, 1984.




\end{thebibliography}
\end{document}